\documentclass[12pt]{article}
\newcommand{\nc}{\newcommand}
\nc{\rnc}{\renewcommand}
\rnc{\baselinestretch}{1.1}
\rnc{\arraystretch}{0.9}
\rnc{\thesubsubsection}{\arabic{subsubsection}.}
\usepackage{pstricks}
\usepackage{latexsym}
\usepackage{graphicx}
\newtheorem{theorem}{Theorem}
\rnc{\thetheorem}{\arabic{theorem}.}
\rnc{\labelenumi}{(\roman{enumi})}

\nc{\be}{\begin{equation}}
\nc{\ee}{\end{equation}}
\nc{\bc}{\begin{center}}
\nc{\ec}{\end{center}}
\nc{\bpic}{\begin{picture}}
\nc{\epic}{\end{picture}}
\nc{\ba}[1]{\begin{array}{@{}#1@{}}}
\nc{\ea}{\end{array}}
\nc{\bea}{\begin{eqnarray}}
\nc{\eea}{\end{eqnarray}}
\nc{\ds}{\displaystyle}
\nc{\ts}{\textstyle}
\rnc{\ss}{\scriptstyle}
\nc{\sss}{\scriptscriptstyle}
\nc{\ru}[1]{\rule[-#1ex]{0ex}{#1ex}}
\nc{\vs}[1]{\vspace*{#1mm}}
\nc{\hs}[1]{\hspace*{#1mm}}
\font\tenmsb=msbm10 scaled \magstep1
\font\sevenmsb=msbm7 scaled \magstep1
\font\fivemsb=msbm5 scaled \magstep1
\newfam\msbfam
\textfont\msbfam=\tenmsb
\scriptfont\msbfam=\sevenmsb
\scriptscriptfont\msbfam=\fivemsb
\def\Bbb#1{{\fam\msbfam\relax#1}}
\nc{\p}[2]{\makebox(0,0)[#1]{$#2$}}
\nc{\pp}[2]{\makebox(0,0)[#1]{$\ss#2$}}
\nc{\ppp}[2]{\makebox(0,0)[#1]{$\sss#2$}}
\nc{\text}[6]{\begin{picture}(#1,#2)\put(#3,#4){\p{#5}{\ds#6}}\end{picture}}
\rnc{\a}{\alpha}
\nc{\A}{\mathcal{A}}
\rnc{\b}{\beta}
\nc{\B}{\mathcal{B}}
\nc{\C}{\Bbb C}
\rnc{\d}{\delta}
\nc{\D}{\Delta}
\nc{\E}{\mathcal{E}}
\nc{\e}{\epsilon}
\nc{\g}{\gamma}
\nc{\hb}{\bar{h}}
\nc{\HH}{\mathcal{H}}
\nc{\Hb}{\bar{H}}
\nc{\La}{A_n}
\nc{\Lp}{\mathcal{L}_\mathcal{P}}
\rnc{\l}{\lambda}
\nc{\m}{\mu}
\nc{\N}{\Bbb N}
\rnc{\O}{\Omega}
\nc{\PP}{\mathcal{P}}
\nc{\PO}{\mathcal{P^\circ}}
\rnc{\P}{\Bbb P}
\nc{\R}{\Bbb R}
\nc{\s}{\sigma}
\rnc{\S}{\mathcal{S}}
\rnc{\t}{\!\times\!}
\nc{\V}{\mathcal{V}}
\nc{\VVb}{\,\bar{\!\mathcal{V}}}
\nc{\Vb}{\,\bar{\!V}}
\nc{\vb}{\bar{v}}
\nc{\Z}{\Bbb Z}
\nc{\mi}{\!-\!}
\nc{\pl}{\!+\!}
\nc{\plmi}{\!\pm\!}
\nc{\mipl}{\!\mp\!}
\nc{\OP}{\mathrm{OP}}
\nc{\EM}{\mathrm{EM}}
\nc{\CEM}{\mathrm{CEM}}
\nc{\ASM}{\mathrm{ASM}}
\nc{\ASMO}{\mathrm{ASM}^\circ}
\nc{\CSM}{\mathrm{CSM}}
\nc{\PM}{\mathrm{PM}}
\nc{\SMS}{\mathrm{SMS}}
\nc{\SMSO}{\mathrm{SMS}^\circ}
\nc{\FPL}{\mathrm{FPL}}
\nc{\LP}{\mathrm{LP}}
\nc{\MT}{\mathrm{MT}}
\nc{\aff}{\mathrm{aff}}
\nc{\ver}{\mathrm{vert}}
\nc{\VT}[5]{\bpic(2,2)\put(1,0.5){\line(0,1){1}}\put(0.5,1){\line(1,0){1}}
\multiput(0.5,1)(1,0){2}{\ppp{}{\bullet}}\multiput(1,0.5)(0,1){2}{\ppp{}{\bullet}}
\put(0.35,1){\ppp{r}{#2}}\put(1,0.3){\ppp{t}{#3}}
\put(1.65,1){\ppp{l}{#4}}\put(1,1.7){\ppp{b}{#5}}\put(1,-0.2){\pp{t}{(#1)}}\epic}
\nc{\bline}[1]{\psline[linecolor=blue,linewidth=1.2pt,linearc=#1]}
\nc{\Vi}{\psline[linewidth=0.5pt](0.5,0)(0.5,1)\psline[linewidth=0.5pt](0,0.5)(1,0.5)}
\nc{\Vv}{\bline{0.2}(-0.5,0)(0.5,0)}
\nc{\Vii}{\bline{0.2}(0,-0.5)(0,0.5)}
\nc{\Viv}{\bline{0.15}(-0.5,0)(0,0)(0,0.5)}
\nc{\Viii}{\bline{0.15}(0,-0.5)(0,0)(0.5,0)}
\nc{\Vx}{\bline{0.17}(-0.5,0)(0,0)(0,0.5)\bline{0.17}(0,-0.5)(0,0)(0.5,0)}
\nc{\Vxiv}{\multips(-0.5,-0.05)(0,0.1){2}{\bline{0.2}(0,0)(1,0)}}
\nc{\Vvi}{\multips(-0.05,-0.5)(0.1,0){2}{\bline{0.2}(0,0)(0,1)}}
\nc{\Vxii}{\bline{0.1}(-0.5,0.05)(-0.05,0.05)(-0.05,0.5)\bline{0.15}(-0.5,-0.05)(0.05,-0.05)(0.05,0.5)}
\nc{\Vviii}{\bline{0.1}(0.05,-0.5)(0.05,-0.05)(0.5,-0.05)\bline{0.15}(-0.05,-0.5)(-0.05,0.05)(0.5,0.05)}
\nc{\Vvii}{\bline{0.1}(-0.05,-0.5)(-0.05,-0.1)(0,0.1)(0,0.5)\bline{0.15}(0.05,-0.5)(0.05,0)(0.5,0)}
\nc{\Vix}{\bline{0.1}(0,-0.5)(0,-0.1)(0.05,0.1)(0.05,0.5)\bline{0.15}(-0.5,0)(-0.05,0)(-0.05,0.5)}
\nc{\Vxi}{\bline{0.1}(-0.5,0)(-0.1,0)(0.1,0.05)(0.5,0.05)\bline{0.15}(0,-0.5)(0,-0.05)(0.5,-0.05)}
\nc{\Vxiii}{\bline{0.1}(-0.5,-0.05)(-0.1,-0.05)(0.1,0)(0.5,0)\bline{0.15}(-0.5,0.05)(0,0.05)(0,0.5)}
\nc{\Vxv}{\bline{0.2}(-0.05,-0.5)(-0.05,-0.1)(0.05,0.1)(0.05,0.5)\bline{0.15}(-0.5,0)(-0.05,0)(-0.05,0.5)\bline{0.15}(0.05,-0.5)(0.05,0)(0.5,0)}
\nc{\Vxviii}{\bline{0.2}(-0.5,-0.05)(-0.1,-0.05)(0.1,0.05)(0.5,0.05)\bline{0.15}(-0.5,0.05)(0,0.05)(0,0.5)\bline{0.15}(0,-0.5)(0,-0.05)(0.5,-0.05)}
\nc{\Vxvi}{\bline{0.2}(-0.5,0)(0,0)(0,0.5)\bline{0.1}(0.05,-0.5)(0.05,-0.05)(0.5,-0.05)\bline{0.18}(-0.05,-0.5)(-0.05,0.05)(0.5,0.05)}
\nc{\Vxvii}{\bline{0.2}(0,-0.5)(0,0)(0.5,0)\bline{0.1}(-0.5,0.05)(-0.05,0.05)(-0.05,0.5)\bline{0.18}(-0.5,-0.05)(0.05,-0.05)(0.05,0.5)}
\nc{\Vxix}{\bline{0.2}(0.05,-0.5)(0.05,-0.05)(0.5,-0.05)\bline{0.3}(-0.05,-0.5)(-0.05,0.05)(0.5,0.05)
\bline{0.2}(-0.5,0.05)(-0.05,0.05)(-0.05,0.5)\bline{0.3}(-0.5,-0.05)(0.05,-0.05)(0.05,0.5)}
\nc{\rline}[1]{\psline[linecolor=red,linewidth=1.2pt,linearc=#1]}
\nc{\VVi}{\rline{0.1}(0.5,0.05)(0.05,0.05)(0.05,0.5)\rline{0.15}(0.5,-0.05)(-0.05,-0.05)(-0.05,0.5)}
\nc{\VVxiv}{\rline{0.1}(-0.5,0.05)(-0.05,0.05)(-0.05,0.5)\rline{0.15}(-0.5,-0.05)(0.05,-0.05)(0.05,0.5)}
\nc{\VVvi}{\rline{0.1}(0.05,-0.5)(0.05,-0.05)(0.5,-0.05)\rline{0.15}(-0.05,-0.5)(-0.05,0.05)(0.5,0.05)}
\nc{\VVxix}{\rline{0.1}(-0.05,-0.5)(-0.05,-0.05)(-0.5,-0.05)\rline{0.15}(0.05,-0.5)(0.05,0.05)(-0.5,0.05)}
\nc{\VVviii}{\multips(-0.05,-0.5)(0.1,0){2}{\rline{0.2}(0,0)(0,1)}}
\nc{\VVxii}{\multips(-0.5,-0.05)(0,0.1){2}{\rline{0.2}(0,0)(1,0)}}
\nc{\VVvii}{\rline{0.1}(-0.05,-0.5)(-0.05,-0.1)(0,0.1)(0,0.5)\rline{0.15}(0.05,-0.5)(0.05,0)(0.5,0)}
\nc{\VVix}{\rline{0.1}(-0.5,0)(-0.1,0)(0.1,0.05)(0.5,0.05)\rline{0.15}(0,-0.5)(0,-0.05)(0.5,-0.05)}
\nc{\VVxi}{\rline{0.1}(0,-0.5)(0,-0.1)(0.05,0.1)(0.05,0.5)\rline{0.15}(-0.5,0)(-0.05,0)(-0.05,0.5)}
\nc{\VVxiii}{\rline{0.1}(-0.5,-0.05)(-0.1,-0.05)(0.1,0)(0.5,0)\rline{0.15}(-0.5,0.05)(0,0.05)(0,0.5)}
\nc{\VViii}{\rline{0.1}(-0.05,0.5)(-0.05,0.1)(0,-0.1)(0,-0.5)\rline{0.15}(0.05,0.5)(0.05,0)(0.5,0)}
\nc{\VViv}{\rline{0.1}(-0.5,0)(-0.1,0)(0.1,-0.05)(0.5,-0.05)\rline{0.15}(0,0.5)(0,0.05)(0.5,0.05)}
\nc{\VVxvi}{\rline{0.1}(0,0.5)(0,0.1)(0.05,-0.1)(0.05,-0.5)\rline{0.15}(-0.5,0)(-0.05,0)(-0.05,-0.5)}
\nc{\VVxvii}{\rline{0.1}(-0.5,0.05)(-0.1,0.05)(0.1,0)(0.5,0)\rline{0.15}(-0.5,-0.05)(0,-0.05)(0,-0.5)}
\nc{\VVii}{\rline{0.15}(0,-0.5)(0,-0.05)(0.5,-0.05)\rline{0.15}(0,0.5)(0,0.05)(0.5,0.05)}
\nc{\VVv}{\rline{0.15}(-0.5,0)(-0.05,0)(-0.05,0.5)\rline{0.15}(0.05,0.5)(0.05,0)(0.5,0)}
\nc{\VVxv}{\rline{0.15}(0.05,-0.5)(0.05,0)(0.5,0)\rline{0.15}(-0.5,0)(-0.05,0)(-0.05,-0.5)}
\nc{\VVxviii}{\rline{0.15}(-0.5,0.05)(0,0.05)(0,0.5)\rline{0.15}(-0.5,-0.05)(0,-0.05)(0,-0.5)}
\nc{\VVxa}{\rline{0.17}(-0.5,0)(0,0)(0,0.5)\rline{0.17}(0,-0.5)(0,0)(0.5,0)}
\nc{\VVxb}{\rline{0.17}(-0.5,0)(0,0)(0,-0.5)\rline{0.17}(0,0.5)(0,0)(0.5,0)}

\begin{document}
\pagestyle{empty}
\bc
\vs{5}
\textbf{\large Higher Spin Alternating Sign Matrices}\\
\bigskip\bigskip
Roger E. Behrend and Vincent A. Knight\\
\bigskip
\textit{School of Mathematics, Cardiff University,\\Cardiff, CF24 4AG, UK}\\
\smallskip
{\footnotesize\tt behrendr@cardiff.ac.uk, knightva@cardiff.ac.uk}\\[28mm]
\begin{abstract}
\noindent We define a higher spin alternating sign matrix to be an integer-entry
square matrix in which, for a nonnegative integer $r$, all complete row and column sums are $r$,
and all partial row and column sums extending from each end of the row or column are nonnegative.
Such matrices correspond to configurations of spin $r/2$ statistical mechanical vertex models with
domain-wall boundary conditions.
The case $r=1$ gives standard alternating sign matrices, while the
case in which all matrix entries are nonnegative gives semimagic
squares.  We show that the higher spin alternating sign matrices of
size~$n$ are the integer points of the $r$-th dilate of an integral convex polytope of dimension $(n\mi1)^2$ whose vertices are
the standard alternating sign matrices of size~$n$.
It then follows that, for fixed $n$, these matrices are enumerated by an
Ehrhart polynomial in $r$.\\[5mm]
\noindent{\footnotesize\emph{Keywords:} alternating sign matrix, semimagic square, convex polytope,
higher spin vertex model}
\end{abstract}\ec
\newpage
\pagestyle{plain}
\subsubsection{Introduction}
Alternating sign matrices are mathematical objects with intriguing combinatorial properties and
important connections to mathematical physics, and the primary aim of this paper is
to introduce natural generalizations of these matrices which also seem to
display interesting such properties and connections.

Alternating sign matrices were first defined in~\cite{MilRobRum82},
and the significance of their connection with mathematical physics first became apparent in~\cite{Kup96},
in which a determinant formula for the partition function of an integrable statistical mechanical
model, and a simple correspondence between configurations of that model and alternating sign matrices, were
used to prove the validity of a previously-conjectured enumeration formula.
For reviews of this and related areas, see for example~\cite{Bre99,BrePro99,Rob91,Zei05}.
Such connections with integrable statistical mechanical models have since been used extensively to derive formulae for
further cases of refined, weighted or symmetry-class enumeration of alternating sign matrices,
as done for example in~\cite{ColPro05,Kup02,RazStr04,Zei96b}.

The statistical mechanical model used in all of these cases is the six-vertex model (with certain boundary conditions),
which is intrinsically related to the spin~$1/2$, or two dimensional, irreducible representation of the Lie algebra $sl(2,\C)$.
For a review of this area, see for example~\cite{GomRuiSie96}.
In this paper, we consider configurations of statistical mechanical vertex models (again with certain boundary conditions)
related to the spin~$r/2$ representation of $sl(2,\C)$, for all nonnegative integers~$r$, these being in simple correspondence
with matrices which we term higher spin alternating sign matrices.  Determinant formulae for the partition
functions of these models have already been obtained in~\cite{CarFodKit06}, thus for example answering Question~22 of~\cite{Kup02}
on whether such formulae exist.

Although we were originally motivated to consider higher spin alternating sign matrices through this
connection with statistical mechanical lattice models, these matrices are natural generalizations
of standard alternating sign matrices in their own right, and appear to have important combinatorial properties.
Furthermore, they generalize not only standard alternating sign matrices, but also other
much-studied combinatorial objects, namely semimagic squares.

Semimagic squares are simply nonnegative integer-entry
square matrices in which all complete row and column sums are equal.
They are thus the integer points of the integer dilates of the convex polytope of
nonnegative real-entry, fixed-size square matrices in which all complete row and column sums are~1,
a fact which leads to enumerative results for the case of fixed
size.  For reviews of this area, see for example~\cite[Ch.~6]{BecRob07} or~\cite[Sec.~4.6]{Sta86}.
In this paper, we introduce an analogous convex polytope,
which was independently defined and studied in~\cite{Str07}, and
for which the integer points of the integer dilates
are the higher spin alternating sign matrices of fixed size.

We define higher spin alternating sign matrices in Section~2, after which this paper then divides into two
essentially independent parts: Sections 3, 4 and 5, and Sections 6, 7 and~8.  In Sections 3, 4 and 5,
we define and discuss various combinatorial objects which are in bijection with higher spin alternating sign matrices,
and which generalize previously-studied objects which are in bijection with standard alternating sign matrices.
In Sections 6, 7 and 8, we define and study the convex polytope which is related to higher spin alternating sign matrices,
and we obtain certain enumerative formulae for the case of fixed size.  We then end the paper in Section~9 with a discussion
of possible further research.

Finally in this introduction, we note that standard alternating sign matrices are related to many further fascinating
results and conjectures
in combinatorics and mathematical physics beyond those already mentioned or directly relevant to this paper.
For example, in combinatorics
it is known that the numbers of standard alternating sign matrices, descending plane partitions,
and totally symmetric self-complementary plane partitions of certain sizes are all equal, but
no bijective proofs of these equalities have yet been found.
Moreover, further equalities between the cardinalities of
certain subsets of these three objects have been conjectured, some over two decades ago, and many of these remain unproved.
See for example~\cite{And79,And94,Dif06,Dif07,Ish06a,Ish06b,MilRobRum83,MilRobRum86}.
Meanwhile, in mathematical physics, extensive work has been done recently on so-called Razumov-Stroganov-type
results and conjectures.  These give surprising equalities between
numbers of certain alternating sign matrices or plane partitions, and
entries of eigenvectors related to certain statistical mechanical models.  See for example~\cite{Deg05,Deg07} and references therein.

\emph{Notation.} \ Throughout this paper, $\P$ denotes the set of positive integers, $\N$ denotes the set of nonnegative integers,
$[m,n]$ denotes the set $\{m,m\pl1,\ldots,n\}$ for any $m,n\in\Z$, with $[m,n]=\emptyset$ for $n<m$, and $[n]$ denotes the set
$[1,n]$ for any $n\in\Z$. The notation $\R_{(0,1)}$ and $\R_{[0,1]}$ will be used for the
open and closed intervals of real numbers between $0$ and $1$.
For a finite set~$T$, $|T|$ denotes the cardinality of~$T$.

\subsubsection{Higher Spin Alternating Sign Matrices}
In this section, we define higher spin alternating sign matrices, describe some of their basic properties,
give an enumeration table, and introduce an example.

For $n\in\P$ and $r\in\N$, let the set of \emph{higher spin alternating sign matrices} of size $n$
with line sum $r$ be
\begin{equation}\label{ASM}\ba{l}\!\!\ASM(n,r)\,:=\\[2mm]
\hs{2.8}\left\{\!A\!=\!\!\left(\ba{c@{\;}c@{\;}c}
A_{11}&\ldots&A_{1n}\\
\vdots&&\vdots\\
A_{n1}&\ldots&A_{nn}\ea\right)\!\in\Z^{n\times n}\:\left|\;\,\ba{l}
\bullet\ \sum_{j'=1}^n\!A_{ij'}=\sum_{i'=1}^n\!A_{i'\!j}=r\mbox{ \ for all }i,j\in[n]\!\\[2.2mm]
\bullet\ \sum_{j'=1}^jA_{ij'\!}\ge0\mbox{ \ for all }i\in[n],\;j\in[n\mi1]\\[2.2mm]
\bullet\ \sum_{j'=j}^nA_{ij'\!}\ge0\mbox{ \ for all }i\in[n],\;j\in[2,n]\\[2.2mm]
\bullet\ \sum_{i'=1}^iA_{i'\!j}\ge0\mbox{ \ for all }i\in[n\mi1],\;j\in[n]\\[2.2mm]
\bullet\ \sum_{i'=i}^nA_{i'\!j}\ge0\mbox{ \ for all }i\in[2,n],\;j\in[n]\ea\right.\right\}\!.\ea\!\!\!\!\end{equation}

In other words, $\ASM(n,r)$ is the set of $n\t n$ integer-entry matrices
for which all complete row and column sums are $r$,
and all partial row and column sums extending from each end of the row or column are nonnegative.
As will be explained in Section~3, a line sum of $r$ corresponds to a spin of $r/2$.
The set $\ASM(n,r)$ can also be written as
\begin{equation}\ba{l}\!\!\ASM(n,r)\,=\\[2mm]
\hs{1.2}\left\{\!A\!=\!\!\left(\ba{c@{\;}c@{\;}c}
A_{11}&\ldots&A_{1n}\\
\vdots&&\vdots\\
A_{n1}&\ldots&A_{nn}\ea\right)\!\in\Z^{n\times n}\:\left|\;\ba{l}
\bullet\ \sum_{j'=1}^nA_{ij'}=\sum_{i'=1}^nA_{i'\!j}=r\mbox{ \ for all }i,j\in[n]\\[2.2mm]
\bullet\ 0\le\sum_{j'=1}^jA_{ij'\!}\le r\mbox{ \ for all }i\in[n],\;j\in[n\mi1]\!\\[2.2mm]
\bullet\ 0\le\sum_{i'=1}^iA_{i'\!j}\le r\mbox{ \ for all }i\in[n\mi1],\;j\in[n]\ea\right.\right\}\!.\ea\!\!\!\!\end{equation}

It follows that each entry of any matrix of $\ASM(n,r)$ is between $-r$ and~$r$,
and that if the entry is in the first or last row or column, then it is
between $0$ and $r$.

A running example will be the matrix
\begin{equation}\label{ExA}A\:=\left(\mbox{\footnotesize$
\ba{ccccc}0&1&1&0&0\\
1&-1&0&2&0\\
0&1&1&-2&2\\
1&0&0&1&0\\
0&1&0&1&0\ea$}\right)\,\in\,\ASM(5,2).\end{equation}

Defining
\begin{equation}\label{SMS}\SMS(n,r):=\{A\in\ASM(n,r)\mid A_{ij}\ge0\mbox{ \ for each }i,j\in[n]\},\end{equation}
it can be seen that this is the set of \emph{semimagic squares} of size $n$ with line sum $r$, i.e., nonnegative integer-entry
$n\t n$ matrices in which all complete row and column sums are~$r$. For example, $\SMS(n,1)$ is the set of
$n\t n$ permutation matrices, so that
\begin{equation}\label{SMS1}|\SMS(n,1)|\,=\,n!\end{equation}
Early studies of semimagic squares appear
in~\cite{AnaDumGup66,Mac15}. For further information and references,
see for example~\cite[Ch.~6]{BecRob07}, \cite{Ehr73}, \cite{Spe80}, \cite{Sta73}, \cite[Sec.~4.6]{Sta86} and~\cite[Sec.~5.5]{Sta99}.

It can also be seen that $\ASM(n,1)$ is the set of \emph{standard alternating sign matrices}
of size $n$, i.e.,
$n\times n$ matrices in which each entry is $0$, $1$ or~$-1$, each row and column contains at least one nonzero entry,
and along each row and column the nonzero entries alternate in sign, starting and finishing with a~1.
Standard alternating sign matrices were first defined and studied in~\cite{MilRobRum82,MilRobRum83}. For further
information, connections to related subjects, and references see for example~\cite{Bre99,BrePro99,Deg07,Pro01,Rob91,Zei05}.

We refer to $\ASM(n,r)$ as a set of
`higher spin alternating sign matrices' for any $n\in\P$ and $r\in\N$, although we realize that this could be slightly misleading since
the `alternating sign' property applies only to the standard case $r=1$, and the spin~$r/2$ is only `higher' for
cases with $r\ge2$.  Nevertheless, we still feel that this is the most natural choice of terminology.

Some cardinalities of $\ASM(n,r)$, many of them computer-generated, are shown in Table~\ref{ASMnr}.

\begin{table}[h]\centering
$\ba{@{\;}r@{\;\;\;}|@{\;\;}r@{\;\;\;\;\;\;}c@{\;\;\;\;}c@{\;\;\;\;}c@{\;\;\;\;}c@{\;}}
\rule{0ex}{2.4ex}&r\!=\!0\!\!\!&1&2&3&4\\[1.5mm]
\hline
\rule{0ex}{3.5ex}n\!=\!1&1&1&1&1&1\\[1.5mm]
2&1&2&3&4&5\\[1.5mm]
3&1&7&26&70&155\\[1.5mm]
4&1&42&628&5102&28005\\[1.5mm]
5&1&429&41784&1507128&28226084\\[1.5mm]
6&1&7436&7517457&1749710096&152363972022
\ea$
\caption{\protect\rule{0ex}{2.5ex}$|\ASM(n,r)|$ \ for $n\in[6]$, $r\in[0,4]$.\label{ASMnr}}
\end{table}

Apart from the trivial formulae $|\ASM(0,n)|=1$ (since $\ASM(0,n)$ contains only the $n\times n$ zero matrix),
$|\ASM(1,r)|=1$ (since $\ASM(1,r)=\{(r)\}$), and $|\ASM(2,r)|=r\pl1$
(since $\ASM(2,r)=\Bigl\{\Bigl(\ba{cc}i&r\mi i\\[-1mm]r\mi i&i\ea\Bigr)\Bigm|i\in[0,r]\Bigr\}=\SMS(2,r)$),
the only previously-known formula for a special case of $|\ASM(n,r)|$ is
\begin{equation}\label{ASM1}|\ASM(n,1)|\;=\;\:\prod_{i=0}^{n\mi1}\!\frac{(3i\pl1)!}{(n\pl i)!}\,,\end{equation}
for standard alternating sign matrices with any $n\in\P$.  This formula was conjectured in~\cite{MilRobRum82,MilRobRum83},
and eventually proved, using different methods, in~\cite{Zei96a} and~\cite{Kup96}.  It has also
been proved using a further method in~\cite{Fis07}, and, using a method related to that of~\cite{Kup96},
in~\cite{ColPro05}.

\subsubsection{Edge Matrix Pairs and Higher Spin Vertex Model Configurations}
In this section, we show that there is a simple bijection between higher spin
alternating sign matrices and configurations of higher spin statistical mechanical vertex models with
domain-wall boundary conditions, and we discuss some properties of these vertex models.

For $n\in\P$ and $r\in\N$, define the set of \emph{edge matrix pairs} as
\begin{equation}\label{EM}\ba{l}\EM(n,r)\,:=\\[2mm]
\;\Biggl\{(H,V)\!=\!\left(\!
\left(\ba{c@{\;}c@{\;}c}H_{10}&\ldots&H_{1n}\\\vdots&&\vdots\\H_{n0}&\ldots&H_{nn}\ea\right)\!,\!
\left(\ba{c@{\;}c@{\;}c}V_{01}&\ldots&V_{0n}\\\vdots&&\vdots\\V_{n1}&\ldots&V_{nn}\ea\right)\!\right)\!
\in[0,r]^{n\times(n\pl1)}\times[0,r]^{(n\pl1)\times n}\:\Bigg|\\[8mm]
\qquad H_{i0}=V_{0j}=0,\ H_{in}=V_{nj}=r,\ H_{i,j\mi1}\pl V_{ij}=V_{i\mi1,j}\pl H_{ij},
\mbox{ \ for all }i,j\in[n]\Biggr\}\!.\!\!\!\ea\end{equation}
We shall refer to $H$ as a \emph{horizontal edge matrix} and $V$ as a \emph{vertical edge matrix}.
It can be checked that there is a bijection between $\ASM(n,r)$ and $\EM(n,r)$
in which the edge matrix pair $(H,V)$ which corresponds to the higher spin alternating sign matrix $A$ is given by
\begin{equation}\label{ASMToHV}\ba{c}\ds H_{ij}=\sum_{j'=1}^jA_{ij'}\,,\mbox{ \ for each }i\in[n],\ j\in[0,n]\\[5mm]
\ds V_{ij}=\sum_{i'=1}^iA_{i'j}\,,\mbox{ \ for each }i\in[0,n],\ j\in[n],\ea\end{equation}
and inversely,
\begin{equation}\label{HVToASM}A_{ij}\,=\,H_{ij}-H_{i,j\mi1}\,=\,V_{ij}-V_{i\mi1,j}\,,\mbox{ \ for each }i,j\in[n].\end{equation}
Thus, $H$ is the \emph{column sum matrix} and $V$ is the \emph{row sum matrix} of $A$.
The correspondence between standard alternating sign matrices and edge matrix pairs was first
identified in~\cite{RobRum86}.

It can be seen that for each $(H,V)\in\EM(n,r)$ and $i,j\in[0,n]$,
$\sum_{i'=1}^nH_{i'j}=jr$ and $\sum_{j'=1}^nV_{ij'}=ir$, so that
\begin{equation}\label{HVSum}\sum_{i,j=1}^nH_{ij}\,=\,\sum_{i,j=1}^nV_{ij}\,=\,n(n\pl1)r/2\,.\end{equation}

The edge matrix pair which corresponds to the running example~(\ref{ExA}) is
\begin{equation}\label{ExHV}(H,V)\,=\left(\left(\mbox{\footnotesize$\ba{cccccc}0&0&1&2&2&2\\
0&1&0&0&2&2\\
0&0&1&2&0&2\\
0&1&1&1&2&2\\
0&0&1&1&2&2\ea$}\right),
\left(\mbox{\footnotesize$\ba{ccccc}0&0&0&0&0\\
0&1&1&0&0\\
1&0&1&2&0\\
1&1&2&0&2\\
2&1&2&1&2\\
2&2&2&2&2\ea$}\right)\right).\end{equation}

\begin{figure}[t]
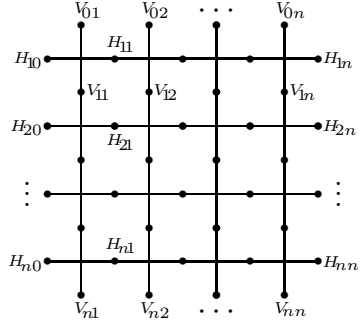
\centering
\bpic(5,4.6)\multiput(0.5,1)(0,1){4}{\line(1,0){4}}\multiput(1,0.5)(1,0){4}{\line(0,1){4}}
\multiput(0.5,1)(1,0){5}{\ppp{}{\bullet}}\multiput(0.5,2)(1,0){5}{\ppp{}{\bullet}}\multiput(0.5,3)(1,0){5}{\ppp{}{\bullet}}
\multiput(0.5,3)(1,0){5}{\ppp{}{\bullet}}\multiput(0.5,4)(1,0){5}{\ppp{}{\bullet}}
\multiput(1,0.5)(0,1){5}{\ppp{}{\bullet}}\multiput(2,0.5)(0,1){5}{\ppp{}{\bullet}}\multiput(3,0.5)(0,1){5}{\ppp{}{\bullet}}
\multiput(3,0.5)(0,1){5}{\ppp{}{\bullet}}\multiput(4,0.5)(0,1){5}{\ppp{}{\bullet}}
\put(1.1,4.6){\ppp{b}{V_{\!01}}}\put(2.1,4.6){\ppp{b}{V_{\!02}}}\put(4.1,4.6){\ppp{b}{V_{\!0n}}}
\put(1.07,3.5){\ppp{l}{V_{\!1\!1}}}\put(2.07,3.5){\ppp{l}{V_{\!1\!2}}}\put(4.07,3.5){\ppp{l}{V_{\!1\!n}}}
\put(1.1,0.4){\ppp{t}{V_{\!n\!1}}}\put(2.1,0.4){\ppp{t}{V_{\!n2}}}\put(4.1,0.4){\ppp{t}{V_{\!nn}}}
\put(0.44,3.98){\ppp{r}{H_{\!1\!0}}}\put(0.44,2.98){\ppp{r}{H_{\!20}}}\put(0.44,0.98){\ppp{r}{H_{\!n0}}}
\put(1.59,4.1){\ppp{b}{H_{\!1\!1}}}\put(1.59,2.88){\ppp{t}{H_{\!2\!1}}}\put(1.59,1.12){\ppp{b}{H_{\!n\!1}}}
\put(4.56,3.98){\ppp{l}{H_{\!1\!n}}}\put(4.56,2.98){\ppp{l}{H_{\!2n}}}\put(4.56,0.98){\ppp{l}{H_{\!nn}}}
\put(0.2,2.12){\p{}{\vdots}}\put(4.8,2.12){\p{}{\vdots}}
\put(2.99,0.25){\p{}{\cdots}}\put(2.99,4.7){\p{}{\cdots}}
\epic\vs{-3}
\caption{Assignment of edge matrix entries to lattice edges.\label{VMC}}\end{figure}

\begin{figure}[t]
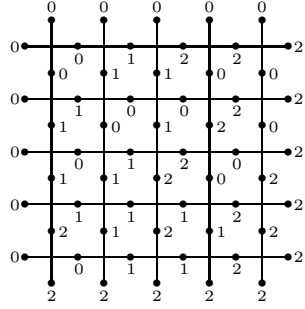
\centering
\bpic(6,6.5)\multiput(0.5,1)(0,1){5}{\line(1,0){5}}\multiput(1,0.5)(1,0){5}{\line(0,1){5}}
\multiput(0.5,1)(1,0){6}{\ppp{}{\bullet}}\multiput(0.5,2)(1,0){6}{\ppp{}{\bullet}}\multiput(0.5,3)(1,0){6}{\ppp{}{\bullet}}
\multiput(0.5,4)(1,0){6}{\ppp{}{\bullet}}\multiput(0.5,5)(1,0){6}{\ppp{}{\bullet}}
\multiput(1,0.5)(0,1){6}{\ppp{}{\bullet}}\multiput(2,0.5)(0,1){6}{\ppp{}{\bullet}}\multiput(3,0.5)(0,1){6}{\ppp{}{\bullet}}
\multiput(4,0.5)(0,1){6}{\ppp{}{\bullet}}\multiput(5,0.5)(0,1){6}{\ppp{}{\bullet}}
\multiput(1,5.65)(1,0){5}{\ppp{b}{0}}\multiput(1,0.35)(1,0){5}{\ppp{t}{2}}\multiput(0.4,1)(0,1){5}{\ppp{r}{0}}\multiput(5.6,1)(0,1){5}{\ppp{l}{2}}
\put(1.13,4.48){\ppp{l}{0}}\put(2.13,4.48){\ppp{l}{1}}\put(3.13,4.48){\ppp{l}{1}}\put(4.13,4.48){\ppp{l}{0}}\put(5.13,4.48){\ppp{l}{0}}
\put(1.13,3.48){\ppp{l}{1}}\put(2.13,3.48){\ppp{l}{0}}\put(3.13,3.48){\ppp{l}{1}}\put(4.13,3.48){\ppp{l}{2}}\put(5.13,3.48){\ppp{l}{0}}
\put(1.13,2.48){\ppp{l}{1}}\put(2.13,2.48){\ppp{l}{1}}\put(3.13,2.48){\ppp{l}{2}}\put(4.13,2.48){\ppp{l}{0}}\put(5.13,2.48){\ppp{l}{2}}
\put(1.13,1.48){\ppp{l}{2}}\put(2.13,1.48){\ppp{l}{1}}\put(3.13,1.48){\ppp{l}{2}}\put(4.13,1.48){\ppp{l}{1}}\put(5.13,1.48){\ppp{l}{2}}
\put(1.52,4.87){\ppp{t}{0}}\put(1.52,3.87){\ppp{t}{1}}\put(1.52,2.87){\ppp{t}{0}}\put(1.52,1.87){\ppp{t}{1}}\put(1.52,0.87){\ppp{t}{0}}
\put(2.52,4.87){\ppp{t}{1}}\put(2.52,3.87){\ppp{t}{0}}\put(2.52,2.87){\ppp{t}{1}}\put(2.52,1.87){\ppp{t}{1}}\put(2.52,0.87){\ppp{t}{1}}
\put(3.52,4.87){\ppp{t}{2}}\put(3.52,3.87){\ppp{t}{0}}\put(3.52,2.87){\ppp{t}{2}}\put(3.52,1.87){\ppp{t}{1}}\put(3.52,0.87){\ppp{t}{1}}
\put(4.52,4.87){\ppp{t}{2}}\put(4.52,3.87){\ppp{t}{2}}\put(4.52,2.87){\ppp{t}{0}}\put(4.52,1.87){\ppp{t}{2}}\put(4.52,0.87){\ppp{t}{2}}
\epic\vs{-3}
\caption{\protect\ru{3}Vertex model configuration for the running example.\label{ExVMC}}\end{figure}

A configuration of a spin $r/2$ statistical mechanical vertex model on an $n\t n$ square
with domain-wall boundary conditions is the assignment, for any $(H,V)\in\EM(n,r)$,
of the horizontal edge matrix entry~$H_{ij}$ to the horizontal edge between lattice points $(i,j)$ and $(i,j\pl1)$, for each $i\in[n]$, $j\in[0,n]$,
and the vertical edge matrix entry~$V_{ij}$ to the vertical edge between lattice points $(i,j)$ and $(i\pl1,j)$, for each $i\in[0,n]$, $j\in[n]$.
Throughout this paper, we use the conventions that
the rows and columns of the lattice are numbered in increasing order from top to bottom, and from left to right,
and that $(i,j)$ denotes the point in row~$i$ and column $j$, i.e., we use matrix-type labeling of lattice points.
The assignment of edge matrix entries to lattice edges is shown diagrammatically in Figure~\ref{VMC}, and
the vertex model configuration for the example of~(\ref{ExHV}) is shown in Figure~\ref{ExVMC}.
The term \emph{domain-wall boundary conditions} refers to the
assignment of $0$ to each edge on the left and upper boundaries of the square,
and of $r$ to each edge on the lower and right boundaries of the square, i.e.,
to the conditions $H_{i0}=V_{0j}=0$ and $H_{in}=V_{nj}=r$ of~(\ref{EM}).
The correspondence between standard alternating sign matrices and configurations of a vertex model with domain-wall boundary conditions
was first identified in~\cite{ElkKupLarPro92}.

We note that in depicting vertex model configurations, it is often standard for
certain numbers of directed arrows, rather than integers in $[0,r]$, to be assigned to lattice edges.
For example, for the case $r=1$, a configuration could be depicted by assigning a
leftward or rightward arrow to the horizontal edge from $(i,j)$ to $(i,j\pl1)$ for $H_{ij}=0$ or
$H_{ij}=1$ respectively,
and assigning a downward or upward arrow to the vertical edge between $(i,j)$ and $(i\pl1,j)$ for
$V_{ij}=0$  or $V_{ij}=1$ respectively.  The condition $H_{i,j\mi1}\pl V_{ij}=V_{i\mi1,j}\pl H_{ij}$ of~(\ref{EM})
then corresponds to arrow conservation at each lattice point
(i.e., that the numbers of arrows into and out of each point are equal), while the domain-wall boundary
conditions correspond to the fact that all arrows on the horizontal or vertical boundaries of the square
point inwards or outwards respectively.

\setlength{\unitlength}{9mm}

\setlength{\unitlength}{7mm}

It is also convenient to define the set of \emph{vertex types}, for a spin $r/2$ statistical mechanical model, as
\begin{equation}\label{Vr}\V(r):=\{(h,v,h',v')\in[0,r]^4\mid h\pl v=h'\pl v'\}.\end{equation}
A vertex type $(h,v,h',v')$ is depicted as \hspace{-5mm}
\setlength{\unitlength}{8mm}
\raisebox{-0.88\unitlength}[1\unitlength][1\unitlength]{
\bpic(2,2)\put(1,0.5){\line(0,1){1}}\put(0.5,1){\line(1,0){1}}
\multiput(0.5,1)(1,0){2}{\ppp{}{\bullet}}\multiput(1,0.5)(0,1){2}{\ppp{}{\bullet}}
\put(0.39,1){\pp{r}{h}}\put(1.63,1){\pp{l}{h'}}
\put(1,0.37){\pp{t}{v}}\put(1.05,1.65){\pp{b}{v'}}\epic},
and it can be seen that for the vertex model configuration associated with $(H,V)\in\EM(n,r)$, the lattice point $(i,j)$
is associated with the vertex type $(H_{i,j\mi1},V_{ij},H_{ij},V_{i\mi1,j})\in\V(r)$, for each $i,j\in[n]$.

The vertex types of $\V(2)$ are shown in Figure~\ref{19V}, where $(1)$--$(19)$ will be used as labels.
The vertex types of $\V(1)$ are $(1)$--$(5)$ and $(10)$ of Figure~\ref{19V}.
\setlength{\unitlength}{8mm}
\begin{figure}[t]\centering
\begin{tabular}{c}
\VT{1}{0}{0}{0}{0}\qquad\VT{2}{0}{1}{0}{1}\qquad\VT{3}{0}{1}{1}{0}\qquad\VT{4}{1}{0}{0}{1}\qquad\VT{5}{1}{0}{1}{0}\qquad\VT{6}{0}{2}{0}{2}\\[7.9mm]
\VT{7}{0}{2}{1}{1}\qquad\VT{8}{0}{2}{2}{0}\qquad\VT{9}{1}{1}{0}{2}\qquad\VT{10}{1}{1}{1}{1}\qquad\VT{11}{1}{1}{2}{0}\qquad\VT{12}{2}{0}{0}{2}\\[7.9mm]
\VT{13}{2}{0}{1}{1}\quad\VT{14}{2}{0}{2}{0}\quad\VT{15}{1}{2}{1}{2}\quad\VT{16}{1}{2}{2}{1}\quad\VT{17}{2}{1}{1}{2}\quad\VT{18}{2}{1}{2}{1}
\quad\VT{19}{2}{2}{2}{2}
\end{tabular}
\caption{\protect\rule[-2.8ex]{0ex}{8.3ex}The 19 vertex types of $\V(2)$.\label{19V}}\end{figure}

For any $r\in\N$, $\V(r)$ can be expressed as the disjoint unions
\begin{equation}\label{VT}\ba{r@{\;}c@{\;\:}l}\V(r)&=&\ds\bigcup_{s=0}^{2r}\;
\Bigl\{(h,s\mi h,h',s\mi h')\Bigm|h,h'\in[\max(0,s\mi r),\min(r,s)]\Bigr\}\\[6mm]
&=&\{(h,v,h',h\pl v\mi h')\mid h,v,h'\in[0,r],\ h\le h'\le v\}\;\cup\\[2.5mm]
&&\{(h,v,h\pl v\mi v',v')\mid h,v,v'\in[0,r],\ v<v'<h\}\;\cup\\[2.5mm]
&&\{(h,h'\pl v'\mi h,h',v')\mid h,h',v'\in[0,r],\ h'<h\le v'\}\;\cup\\[2.5mm]
&&\{(h'\pl v'\mi v,v,h',v')\mid v,h',v'\in[0,r],\ v'\le v<h'\}\,,\ea\end{equation}
so that
\begin{equation}|\V(r)|=2\sum_{s=1}^{r}\!s^2+(r\pl1)^2=
\biggl(\ba{c}r\pl1\\[-0.2mm]3\ea\biggr)
+2\biggl(\ba{c}r\pl2\\[-0.2mm]3\ea\biggr)
+\biggl(\ba{c}r\pl3\\[-0.2mm]3\ea\biggr)=(r\pl1)(2r^2\pl4r\pl3)/3.\end{equation}

It can be seen, using~(\ref{SMS}) and~(\ref{HVToASM}), that a spin $r/2$ vertex model configuration corresponds to a
semimagic square with line sum~$r$
if and only if each of its vertex types is in $\V_{\ss\mathrm{S}}(r):=\{(h,v,h',v')\in\V(r)\,|\,h\le h'\mbox{ (and }v'\le v)\}$.
For example, $\V_{\ss\mathrm{S}}(1)$ consists of (1)--(3), (5) and (10) of Figure~\ref{19V}, and
$\V_{\ss\mathrm{S}}(2)$ consists of (1)--(3), (5)--(8), (10), (11), (14)--(16),
(18) and (19) of Figure~\ref{19V}.
\newpage
By imposing the condition $h\le h'$ on the two disjoint unions of~(\ref{VT}), which in the second case leaves just
the first and fourth sets, it follows that \ru{1.6}$|\V_{\ss\mathrm{S}}(r)|=\sum_{s=1}^{r\pl1}\!s^2=
\biggl(\ba{c}r\pl2\\[-0.2mm]3\ea\biggr)
+\biggl(\ba{c}r\pl3\\[-0.2mm]3\ea\biggr)=(r\pl1)(r\pl2)(2r\pl3)/6$.

For a spin $r/2$ statistical mechanical vertex model, a \emph{Boltzmann weight}
\begin{equation}\label{W}W(r,x,h,v,h',v')\in\C\end{equation}
is defined for each $(h,v,h',v')\in\V(r)$. Here, $x$ is
a complex variable, often called the \emph{spectral parameter}.

For such a model on an $n$ by $n$ square with domain-wall boundary conditions, and an $n\t n$ matrix $z$ with entries
$z_{ij}\in\C$ for $i,j\in[n]$, the \emph{partition function} is
\begin{equation}\label{Z}Z(n,r,z)\;:=\sum_{(H,V)\in\,\EM(n,r)}\;\prod_{i,j=1}^n\,W(r,z_{ij},H_{i,j\mi1},V_{ij},H_{ij},V_{i\mi1,j})\,.\end{equation}
Values of $Z(n,r,z)$ therefore give certain weighted enumerations of the higher spin alternating sign matrices of~$\ASM(n,r)$.
It follows that if there exists $u^{\sss\mathrm{A}}_r\in\C$ such that
\begin{equation}\label{ASMCond}W(r,u_r^{\sss\mathrm{A}},h,v,h',v')\,=\,1\qquad\mbox{for each }(h,v,h',v')\in\V(r),\end{equation}
then
\begin{equation}\label{ZASM}Z(n,r,z)|_{\,\mathrm{each}\,z_{ij}=u_r^{\sss\mathrm{A}}}\,=\,|\ASM(n,r)|\,,\end{equation}
and that if there exists $u_r^{\sss\mathrm{S}}\in\C$ such that
\begin{equation}\label{SMSCond}W(r,u_r^{\sss\mathrm{S}},h,v,h',v')\,=\left\{\ba{l}1,\ h\le h'\\
0,\ h>h'\ea\right.\qquad\mbox{for each }(h,v,h',v')\in\V(r),\end{equation}
then
\begin{equation}\label{ZSMS}Z(n,r,z)|_{\,\mathrm{each}\,z_{ij}=u_r^{\sss\mathrm{S}}}\,=\,|\SMS(n,r)|\,.\end{equation}

The Boltzmann weights~(\ref{W}) are usually assumed to satisfy
the \emph{Yang-Baxter equation} and certain other properties.  See for example~\cite[Ch.~8~\&~9]{Bax82} and~\cite[Ch.~1~\&~2]{GomRuiSie96}.
Such a model is then known as integrable, and is
related to the spin $r/2$ representation, i.e., the irreducible representation with highest weight $r$
and dimension $r\pl1$, of the simple Lie algebra $sl(2,\C)$, or its affine counterpart.  See for example~\cite{Fuc92,FucSch97,GomRuiSie96}.
Each value $i\in[0,r]$, as taken by the edge matrix entries, can thus be associated with
an $sl(2,\C)$ weight $2i\mi r$.  In physics contexts, it is also natural to associate
each $i\in[0,r]$ with a spin value $i\mi r/2$.
The model with $r=1$ is known as the \emph{six-vertex} or \emph{square ice} model, and is
related to the \emph{XXZ spin chain} and the defining spin~$1/2$ representation of $sl(2,\C)$.
Furthermore, the Boltzmann weights for each case with $r>1$ can be obtained from the
Boltzmann weights for the $r=1$ case using a procedure known
as \emph{fusion}. See for example~\cite{KulResSkl81}. Boltzmann weights for the $r=2$ case
are also obtained more directly in~\cite{IdzTokAra94,SogAkuAbe83,ZamFat80}.

For the Boltzmann weights of these models, and for any $x\!=\!(x_1,\ldots,x_n)$, $y\!=\!(y_1,\ldots,y_n)\in\C^n$ with each having distinct entries,
it can be shown that
\begin{equation}\label{ZDet}Z(n,r,z)|_{\,\mathrm{each}\,z_{ij}=x_i\mi y_j}\;=\;F(n,r,x,y)\;\det M(n,r,x,y)\,,\end{equation}
where $M(n,r,x,y)$ is an $nr\times nr$ matrix with
entries $M(n,r,x,y)_{(i,k),(j,l)}=\phi(k\mi l,$ $x_i\mi y_j)$ for each $(i,k),(j,l)\in[n]\t[r]$, and
$F$ and $\phi$ are relatively simple, explicitly-known functions.  This determinant formula for the partition function
with $z_{ij}=x_i\mi y_j$ is proved for the case $r=1$ in~\cite{Ize87,IzeCokKor92}, using results of~\cite{Kor82}, and for each
case with $r>1$ in~\cite{CarFodKit06}, using the $r=1$ result and the fusion procedure.
The formula for the $r=1$ case is also proved in~\cite{BogProZvo02}, using a method different from that of~\cite{Ize87,IzeCokKor92},
while that for each case of $r>1$ was obtained independently of~\cite{CarFodKit06}, but using a similar fusion method,
in~\cite{Beh05}.

If any entries of~$x$, or any entries of~$y$, are equal, then $F(n,r,x,y)$ has a singularity, and $\det M(n,r,x,y)=0$.
However, by taking an appropriate limit as the entries become equal, as done in~\cite{IzeCokKor92} for
$r=1$ and~\cite{CarFodKit06} for $r>1$, a valid alternative formula involving the determinant of an $nr\times nr$ matrix whose
entries are derivatives of the function $\phi$ can be obtained.
For the completely homogeneous case in which all entries of $x$ are equal,
and all entries of~$y$ are equal, with a difference $u$ between the entries of $x$ and $y$,
this matrix has entries $\frac{d^{i\pl j\mi2}}{\rule{0ex}{1.5ex}du^{i\pl j\mi2}}\,\phi(k\mi l,u)$ for each $(i,k),(j,l)\in[n]\t[r]$.

For the case $r=1$, there exists $u^{\sss\mathrm{A}}_1$ such that the Boltzmann weights which lead to a
determinant formula~(\ref{ZDet}) satisfy~(\ref{ASMCond}), so that~(\ref{ZASM}) can be applied.  This is
done in~\cite{Kup96} and~\cite{ColPro05} in order to prove~(\ref{ASM1}).  In~\cite{Kup96}, a choice of $x$ and $y$
which depend on a parameter $\e$ is used, in which $x$ and $y$
each have distinct entries for $\e\ne0$, and $x_i\mi y_j=u^{\sss\mathrm{A}}_1$ for $\e=0$ and each $i,j\in[n]$.  The formula~(\ref{ZDet}) is
then applied with $\e\ne0$, the resulting determinant is evaluated as a product form, and finally the limit $\e\rightarrow0$ is
taken, giving the RHS of~(\ref{ASM1}).  In~\cite{ColPro05}, a
determinant formula for the completely homogeneous case is applied at the outset, and the relation between Hankel determinants
and orthogonal polynomials, together with known properties
of the Continuous Hahn orthogonal polynomials, are then used to evaluate the resulting determinant,
giving the RHS of~(\ref{ASM1}).

For the cases with $r>1$, if there exist values $u^{\sss\mathrm{A}}_r$ and $u^{\sss\mathrm{S}}_r$ such that Boltzmann weights which lead to a
determinant formula~(\ref{ZDet}) satisfy~(\ref{ASMCond}) and~(\ref{SMSCond}), then methods similar to those used in the $r=1$ case could be
applied in an attempt to obtain formulae for $|\ASM(n,r)|$ and $|\SMS(n,r)|$ for fixed~$r$ and variable~$n$.
However, unfortunately it seems that such $u^{\sss\mathrm{A}}_r$ and $u^{\sss\mathrm{S}}_r$ might not exist.

\subsubsection{Lattice Paths}
In this section, we show that there is also a bijection between higher
spin alternating sign matrices and certain sets of lattice paths.

For $n\in\P$ and $r\in\N$, let $\LP(n,r)$ be the set of all sets $P$ of $n r$ directed lattice paths such that
\vs{-5}
\begin{itemize}
\item For each $i\in[n]$, $P$ contains $r$ paths which begin by passing from\\
$(n\pl1,i)$ to $(n,i)$ and end by passing from $(i,n)$ to $(i,n\pl1)$.\vs{-2}
\item Each step of each path of $P$ is either $(-1,0)$ or $(0,1)$.\vs{-2}
\item Different paths of $P$ do not cross.\vs{-2}
\item No more than $r$ paths of $P$ pass along any edge of the lattice.
\end{itemize}

It can be checked that there is a bijection between $\EM(n,r)$ (and hence $\ASM(n,r)$) and
$\LP(n,r)$ in which the edge matrix pair $(H,V)$ which corresponds to the
path set $P$ is given simply by
\begin{equation}\label{PToHV}\ba{r@{\;\;}l}H_{ij}\;=&\mbox{number of paths of }P\mbox{ which pass from }(i,j)\\[1mm]
&\mbox{ to }(i,j\pl1),\mbox{ \ for each }i\in[n],\,j\in[0,n]\\[3mm]
V_{ij}\;=&\mbox{number of paths of }P\mbox{ which pass from }(i\pl1,j)\\[1mm]
&\mbox{ to }(i,j),\mbox{ \ for each }i\in[0,n],\,j\in[n].\ea\end{equation}

\begin{figure}[h]
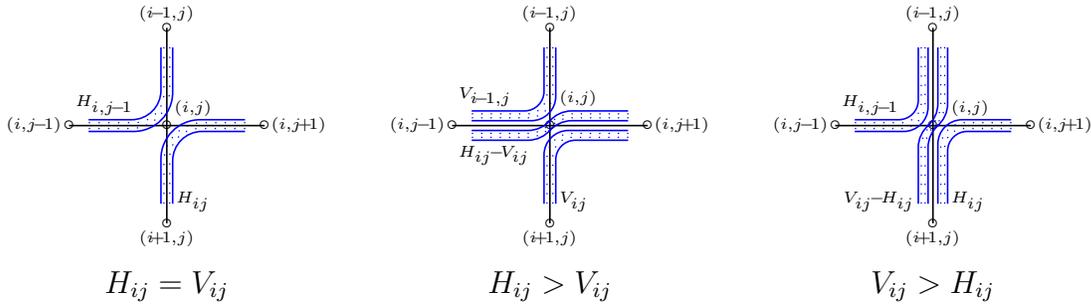
\centering
\pspicture(-1.8,-1.8)(1.8,1.8)
\psline[linewidth=0.6pt](-1,0)(1,0)\psline[linewidth=0.6pt](0,-1)(0,1)
\rput[b](0,-1.8){$H_{ij}=V_{ij}$}
\multirput(-1,0)(1,0){3}{$\ss\circ$}\multirput(0,-1)(0,2){2}{$\ss\circ$}
\rput[r](-1.07,0){$\sss(i,j\mi1)$}\rput[l](1.07,0){$\sss(i,j\pl1)$}
\rput[t](0,-1.07){$\sss(i\pl1,j)$}\rput[b](0,1.07){$\sss(i\mi1,j)$}
\rput[bl](0.08,0.1){$\sss(i,j)$}
\rput[bl](-0.93,0.11){$\sss H_{i,j\mi1}$}\rput[l](0.1,-0.78){$\sss H_{ij}$}
\psline[linearc=.4,linewidth=0.6pt,linecolor=blue](-0.8,-0.06)(0.06,-0.06)(0.06,0.8)
\psline[linearc=.4,linewidth=0.6pt,linecolor=blue,linestyle=dotted](-0.8,-0.03)(0.03,-0.03)(0.03,0.8)
\psline[linearc=.3,linewidth=0.6pt,linecolor=blue,linestyle=dotted](-0.8,0.03)(-0.03,0.03)(-0.03,0.8)
\psline[linearc=.3,linewidth=0.6pt,linecolor=blue](-0.8,0.06)(-0.06,0.06)(-0.06,0.8)
\psline[linearc=.4,linewidth=0.6pt,linecolor=blue](-0.06,-0.8)(-0.06,0.06)(0.8,0.06)
\psline[linearc=.4,linewidth=0.6pt,linecolor=blue,linestyle=dotted](-0.03,-0.8)(-0.03,0.03)(0.8,0.03)
\psline[linearc=.3,linewidth=0.6pt,linecolor=blue,linestyle=dotted](0.03,-0.8)(0.03,-0.03)(0.8,-0.03)
\psline[linearc=.3,linewidth=0.6pt,linecolor=blue](0.06,-0.8)(0.06,-0.06)(0.8,-0.06)
\endpspicture
\quad
\pspicture(-1.8,-1.8)(1.8,1.3)
\psline[linewidth=0.6pt](-1,0)(1,0)\psline[linewidth=0.6pt](0,-1)(0,1)
\rput[b](0,-1.8){$H_{ij}>V_{ij}$}
\multirput(-1,0)(1,0){3}{$\ss\circ$}\multirput(0,-1)(0,2){2}{$\ss\circ$}
\rput[r](-1.07,0){$\sss(i,j\mi1)$}\rput[l](1.07,0){$\sss(i,j\pl1)$}
\rput[t](0,-1.07){$\sss(i\pl1,j)$}\rput[b](0,1.07){$\sss(i\mi1,j)$}
\rput[bl](0.1,0.2){$\sss(i,j)$}
\rput[tl](-0.93,-0.21){$\sss H_{ij}\!-\!V_{ij}$}\rput[bl](-0.93,0.2){$\sss V_{i\mi1,j}$}\rput[l](0.1,-0.78){$\sss V_{ij}$}
\psline[linearc=.25,linewidth=0.6pt,linecolor=blue](-0.8,0.05)(0.06,0.05)(0.06,0.8)
\psline[linearc=.25,linewidth=0.6pt,linecolor=blue,linestyle=dotted](-0.8,0.08)(0.03,0.08)(0.03,0.8)
\psline[linearc=.2,linewidth=0.6pt,linecolor=blue,linestyle=dotted](-0.8,0.12)(-0.03,0.12)(-0.03,0.8)
\psline[linearc=.2,linewidth=0.6pt,linecolor=blue](-0.8,0.15)(-0.06,0.15)(-0.06,0.8)
\psline[linearc=.25,linewidth=0.6pt,linecolor=blue](-0.06,-0.8)(-0.06,-0.05)(0.8,-0.05)
\psline[linearc=.25,linewidth=0.6pt,linecolor=blue,linestyle=dotted](-0.03,-0.8)(-0.03,-0.08)(0.8,-0.08)
\psline[linearc=.2,linewidth=0.6pt,linecolor=blue,linestyle=dotted](0.03,-0.8)(0.03,-0.12)(0.8,-0.12)
\psline[linearc=.2,linewidth=0.6pt,linecolor=blue](0.06,-0.8)(0.06,-0.15)(0.8,-0.15)
\psline[linearc=.3,linewidth=0.6pt,linecolor=blue](-0.8,-0.05)(-0.08,-0.05)(0.08,0.15)(0.8,0.15)
\psline[linearc=.3,linewidth=0.6pt,linecolor=blue,linestyle=dotted](-0.8,-0.08)(-0.08,-0.08)(0.08,0.12)(0.8,0.12)
\psline[linearc=.3,linewidth=0.6pt,linecolor=blue,linestyle=dotted](-0.8,-0.12)(-0.08,-0.12)(0.08,0.08)(0.8,0.08)
\psline[linearc=.3,linewidth=0.6pt,linecolor=blue](-0.8,-0.15)(-0.08,-0.15)(0.08,0.05)(0.8,0.05)
\endpspicture
\quad
\pspicture(-1.8,-1.8)(1.8,1.3)
\psline[linewidth=0.6pt](-1,0)(1,0)\psline[linewidth=0.6pt](0,-1)(0,1)
\rput[b](0,-1.8){$V_{ij}>H_{ij}$}
\multirput(-1,0)(1,0){3}{$\ss\circ$}\multirput(0,-1)(0,2){2}{$\ss\circ$}
\rput[r](-1.07,0){$\sss(i,j\mi1)$}\rput[l](1.07,0){$\sss(i,j\pl1)$}
\rput[t](0,-1.07){$\sss(i\pl1,j)$}\rput[b](0,1.07){$\sss(i\mi1,j)$}
\rput[bl](0.19,0.1){$\sss(i,j)$}
\rput[bl](-0.93,0.11){$\sss H_{i,j\mi1}$}\rput[r](-0.19,-0.78){$\sss V_{ij}\mi H_{ij}$}\rput[l](0.19,-0.78){$\sss H_{ij}$}
\psline[linearc=.25,linewidth=0.6pt,linecolor=blue](-0.8,-0.06)(-0.05,-0.06)(-0.05,0.8)
\psline[linearc=.25,linewidth=0.6pt,linecolor=blue,linestyle=dotted](-0.8,-0.03)(-0.08,-0.03)(-0.08,0.8)
\psline[linearc=.2,linewidth=0.6pt,linecolor=blue,linestyle=dotted](-0.8,0.03)(-0.12,0.03)(-0.12,0.8)
\psline[linearc=.2,linewidth=0.6pt,linecolor=blue](-0.8,0.06)(-0.15,0.06)(-0.15,0.8)
\psline[linearc=.25,linewidth=0.6pt,linecolor=blue](0.05,-0.8)(0.05,0.06)(0.8,0.06)
\psline[linearc=.25,linewidth=0.6pt,linecolor=blue,linestyle=dotted](0.08,-0.8)(0.08,0.03)(0.8,0.03)
\psline[linearc=.2,linewidth=0.6pt,linecolor=blue,linestyle=dotted](0.12,-0.8)(0.12,-0.03)(0.8,-0.03)
\psline[linearc=.2,linewidth=0.6pt,linecolor=blue](0.15,-0.8)(0.15,-0.06)(0.8,-0.06)
\psline[linearc=.3,linewidth=0.6pt,linecolor=blue](-0.05,-0.8)(-0.05,-0.08)(0.15,0.08)(0.15,0.8)
\psline[linearc=.3,linewidth=0.6pt,linecolor=blue,linestyle=dotted](-0.08,-0.8)(-0.08,-0.08)(0.12,0.08)(0.12,0.8)
\psline[linearc=.3,linewidth=0.6pt,linecolor=blue,linestyle=dotted](-0.12,-0.8)(-0.12,-0.08)(0.08,0.08)(0.08,0.8)
\psline[linearc=.3,linewidth=0.6pt,linecolor=blue](-0.15,-0.8)(-0.15,-0.08)(0.05,0.08)(0.05,0.8)
\endpspicture\\[-1mm]
\caption{Path configurations through vertex $(i,j)$ for the cases of (\ref{LinkRules}).\label{PathLink}}\end{figure}

\begin{figure}[t]
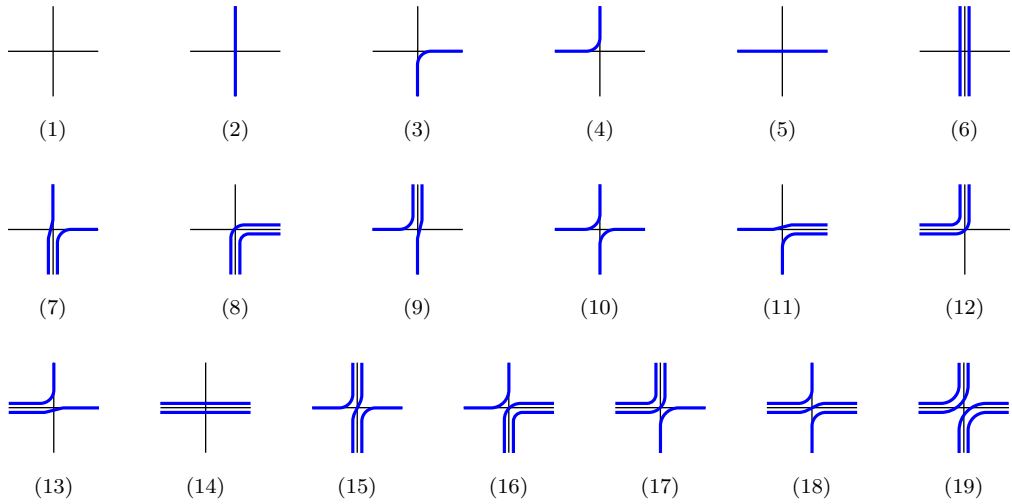
\centering
\begin{tabular}{c}
\pspicture(0,-0.5)(1.51,1.7)\Vi\rput[b](0.5,-0.5){$\ss(1)$}\endpspicture
\pspicture(-0.51,-0.5)(1.51,1)\Vi\rput(0.5,0.5){\Vii}\rput[b](0.5,-0.5){$\ss(2)$}\endpspicture
\pspicture(-0.51,-0.5)(1.51,1)\Vi\rput(0.5,0.5){\Viii}\rput[b](0.5,-0.5){$\ss(3)$}\endpspicture
\pspicture(-0.51,-0.5)(1.51,1)\Vi\rput(0.5,0.5){\Viv}\rput[b](0.5,-0.5){$\ss(4)$}\endpspicture
\pspicture(-0.51,-0.5)(1.51,1)\Vi\rput(0.5,0.5){\Vv}\rput[b](0.5,-0.5){$\ss(5)$}\endpspicture
\pspicture(-0.51,-0.5)(1,1)\Vi\rput(0.5,0.5){\Vvi}\rput[b](0.5,-0.5){$\ss(6)$}\endpspicture\\
\pspicture(0,-0.5)(1.51,1.35)\Vi\rput(0.5,0.5){\Vvii}\rput[b](0.5,-0.5){$\ss(7)$}\endpspicture
\pspicture(-0.51,-0.5)(1.51,1.35)\Vi\rput(0.5,0.5){\Vviii}\rput[b](0.5,-0.5){$\ss(8)$}\endpspicture
\pspicture(-0.51,-0.5)(1.51,1.35)\Vi\rput(0.5,0.5){\Vix}\rput[b](0.5,-0.5){$\ss(9)$}\endpspicture
\pspicture(-0.51,-0.5)(1.51,1.35)\Vi\rput(0.5,0.5){\Vx}\rput[b](0.5,-0.5){$\ss(10)$}\endpspicture
\pspicture(-0.51,-0.5)(1.51,1.35)\Vi\rput(0.5,0.5){\Vxi}\rput[b](0.5,-0.5){$\ss(11)$}\endpspicture
\pspicture(-0.51,-0.5)(1,1.35)\Vi\rput(0.5,0.5){\Vxii}\rput[b](0.5,-0.5){$\ss(12)$}\endpspicture\\
\pspicture(0,-0.6)(1.34,1.35)\Vi\rput(0.5,0.5){\Vxiii}\rput[b](0.5,-0.5){$\ss(13)$}\endpspicture
\pspicture(-0.34,-0.6)(1.34,1.35)\Vi\rput(0.5,0.5){\Vxiv}\rput[b](0.5,-0.5){$\ss(14)$}\endpspicture
\pspicture(-0.34,-0.6)(1.34,1.35)\Vi\rput(0.5,0.5){\Vxv}\rput[b](0.5,-0.5){$\ss(15)$}\endpspicture
\pspicture(-0.34,-0.6)(1.34,1.35)\Vi\rput(0.5,0.5){\Vxvi}\rput[b](0.5,-0.5){$\ss(16)$}\endpspicture
\pspicture(-0.34,-0.6)(1.34,1.35)\Vi\rput(0.5,0.5){\Vxvii}\rput[b](0.5,-0.5){$\ss(17)$}\endpspicture
\pspicture(-0.34,-0.6)(1.34,1.35)\Vi\rput(0.5,0.5){\Vxviii}\rput[b](0.5,-0.5){$\ss(18)$}\endpspicture
\pspicture(-0.34,-0.6)(1,1.35)\Vi\rput(0.5,0.5){\Vxix}\rput[b](0.5,-0.5){$\ss(19)$}\endpspicture\end{tabular}
\caption{Path configurations for the 19 vertex types of $\V(2)$.\label{19VLP}}\end{figure}

\begin{figure}[h]
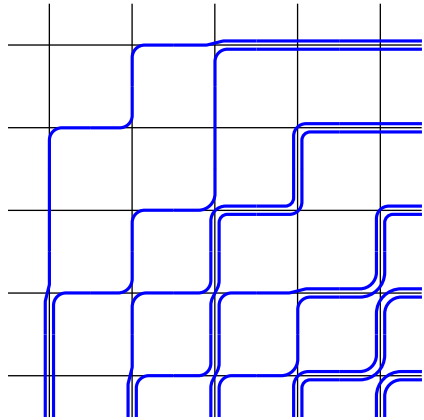
\centering
\pspicture(0.5,0.5)(5.5,5.8)
\multips(0.5,1)(0,1){5}{\psline[linewidth=0.5pt](0,0)(5,0)}\multips(1,0.5)(1,0){5}{\psline[linewidth=0.5pt](0,0)(0,5)}
\rput(2,5){\Viii}\rput(3,5){\Vxi}\rput(4,5){\Vxiv}\rput(5,5){\Vxiv}
\rput(1,4){\Viii}\rput(2,4){\Viv}\rput(3,4){\Vii}\rput(4,4){\Vviii}\rput(5,4){\Vxiv}
\rput(1,3){\Vii}\rput(2,3){\Viii}\rput(3,3){\Vxvi}\rput(4,3){\Vxii}\rput(5,3){\Vviii}
\rput(1,2){\Vvii}\rput(2,2){\Vx}\rput(3,2){\Vxv}\rput(4,2){\Vxi}\rput(5,2){\Vxix}
\rput(1,1){\Vvi}\rput(2,1){\Vvii}\rput(3,1){\Vxv}\rput(4,1){\Vxvi}\rput(5,1){\Vxix}
\endpspicture
\caption{Set of lattice paths for the running example.\label{ExLP}}\end{figure}

For the inverse mapping from $(H,V)$ to $P$,~(\ref{PToHV}) is used to
assign appropriate numbers of path segments to the horizontal and vertical edges of the lattice, and
at each $(i,j)\in[n]\t[n]$, the $H_{i,j\mi1}\pl V_{ij}=V_{i\mi1,j}\pl H_{ij}$
segments on the four neighboring edges are linked without crossing through $(i,j)$ according to the rules that
\begin{equation}\label{LinkRules}
\ba{rl}\bullet&\mbox{If $H_{ij}=V_{ij}$ (and $H_{i,j\mi1}=V_{i\mi1,j}$), then $H_{i,j\mi1}$ paths pass from $(i,j\mi1)$}\\[1mm]
&\mbox{to $(i\mi1,j)$, and $H_{ij}$ paths pass from $(i\pl1,j)$ to $(i,j\pl1)$.}\\[2.5mm]
\bullet&\mbox{If $H_{ij}>V_{ij}$ (and $H_{i,j\mi1}>V_{i\mi1,j}$), then $V_{i\mi1,j}$ paths pass from $(i,j\mi1)$}\\[1mm]
&\mbox{to $(i\mi1,j)$, $H_{ij}\mi V_{ij}=H_{i,j\mi1}\mi V_{i\mi1,j}$ paths pass from $(i,j\mi1)$}\\[1mm]
&\mbox{to $(i,j\pl1)$, and $V_{ij}$ paths pass from $(i\pl1,j)$ to $(i,j\pl1)$.}\\[2.5mm]
\bullet&\mbox{If $V_{ij}>H_{ij}$ (and $V_{i\mi1,j}>H_{i,j\mi1}$), then $H_{i,j\mi1}$ paths pass from $(i,j\mi1)$}\\[1mm]
&\mbox{to $(i\mi1,j)$, $V_{ij}\mi H_{ij}=V_{i\mi1,j}\mi H_{i,j\mi1}$ paths pass from $(i\pl1,j)$}\\[1mm]
&\mbox{to $(i\mi1,j)$, and $H_{ij}$ paths pass from $(i\pl1,j)$ to $(i,j\pl1)$.}\ea
\end{equation}
The three cases of~(\ref{LinkRules}) are shown diagrammatically in Figure~\ref{PathLink},
the path configurations which correspond to the vertex types of $\V(2)$ from Figure~\ref{19V} are shown
in Figure~\ref{19VLP}, and the path set of $\LP(5,2)$ which corresponds to the running example
of~(\ref{ExA}),~(\ref{ExHV}) and Figure~\ref{ExVMC} is shown in Figure~\ref{ExLP}.  In order to
assist in their visualization, some of the path segments
in these diagrams have been shifted slightly away from the lattice edges on which they actually lie.
Also, as indicated in the previous section, we are using matrix-type labeling of lattice points.

\psset{unit=13mm}

\vs{3}
The case $\LP(n,1)$ of path sets for standard alternating sign matrices
is studied in detail in~\cite{Beh07} as a particular case of osculating paths
which start and end at fixed points on the lower and right boundaries of a rectangle.
The correspondence between standard alternating sign matrices and such
osculating paths is also considered in~\cite[Sec.~5]{BouHab95},~\cite[Sec.~2]{Bra97},~\cite[Sec.~9]{EgeRedRya01}
and~\cite[Sec.~IV]{Tam01}.
\newpage
\psset{unit=12mm}

\psset{unit=11mm}

\subsubsection{Further Representations of Higher Spin Alternating Sign Matrices}
In this section, we describe three further combinatorial objects which are in bijection with
higher spin alternating sign matrices: corner sum matrices, monotone triangles and complementary
edge matrix pairs.
These provide generalizations of previously-studied combinatorial objects which are
in bijection with standard alternating sign matrices.

For $n\in\P$ and $r\in\N$, let the set of \emph{corner sum matrices} be
\begin{equation}\label{CSM}\!\!\!\ba{l}\CSM(n,r):=\!
\;\Biggl\{C\!=\!\left(\ba{c@{\;}c@{\;}c}C_{00}&\ldots&C_{0n}\\\vdots&&\vdots\\C_{n0}&\ldots&C_{nn}\ea\right)\!\in
\N^{(n\pl1)\times(n\pl1)}\:\Bigg|\\[9.4mm]
\hs{30}\ba{l}
\bullet\ C_{0k}=C_{k0}=0,\;C_{kn}=C_{nk}=kr,\mbox{ \ for all }k\in[n]\\[2.2mm]
\bullet\ 0\le C_{ij}\mi C_{i,j\mi1}\le r,\;\ 0\le C_{ij}\mi C_{i\mi1,j}\le r,
\mbox{ \ for all }i,j\in[n]\ea\Biggr\}\!.\ea\!\!\!\!\!\!\end{equation}
It can be checked that there is a bijection between $\ASM(n,r)$ and $\CSM(n,r)$
in which the corner sum matrix $C$ which corresponds to the higher spin alternating sign matrix $A$ is given by
\begin{equation}\label{ASMToCSM}\ds C_{ij}\,=\,\sum_{i'=1}^i\sum_{j'=1}^jA_{i'\!j'}\,,\mbox{ \ for each }i,j\in[0,n],\end{equation}
and inversely,
\begin{equation}\label{CSMToASM}A_{ij}\,=\,C_{ij}-C_{i,j\mi1}-C_{i\mi1,j}+C_{i\mi1,j\mi1}\,,\mbox{ \ for each }i,j\in[n].\end{equation}
Combining the bijections~(\ref{ASMToHV},\ref{HVToASM}) between $\EM(n,r)$ and $\ASM(n,r)$,
and~(\ref{ASMToCSM},\ref{CSMToASM}) between $\ASM(n,r)$ and $\CSM(n,r)$,
the corner sum matrix $C$ which corresponds to the edge matrix pair $(H,V)$ is given by
\begin{equation}\label{HVToCSM}\ds C_{ij}\,=\,\sum_{i'=1}^iH_{i'\!j}\,=\,\sum_{j'=1}^jV_{ij'}\,,\mbox{ \ for each }i,j\in[0,n],\end{equation}
and inversely,
\begin{equation}\label{CSMToHV}\ba{r@{}l@{}l}H_{ij}\;&=\,C_{ij}-C_{i\mi1,j}\,,&\mbox{ \ for each }i\in[n],\ j\in[0,n]\\[3mm]
V_{ij}\;&=\,C_{ij}-C_{i,j\mi1}\,,&\mbox{ \ for each }i\in[0,n],\ j\in[n].\ea\end{equation}

The set $\CSM(n,1)$ of corner sum matrices for standard alternating sign matrices
was introduced in~\cite{RobRum86}, and is also considered in~\cite{Pro01}.

The corner sum matrix which corresponds to the running example of~(\ref{ExA}) and~(\ref{ExHV}) is
\begin{equation}\label{ExC}\left(\mbox{\footnotesize$
\ba{cccccc}0&0&0&0&0&0\\
0&0&1&2&2&2\\
0&1&1&2&4&4\\
0&1&2&4&4&6\\
0&2&3&5&6&8\\
0&2&4&6&8&\!10\ea$}\right).\end{equation}

Proceeding now to sets of \emph{monotone triangles},
for $n\in\P$ and $r\in\N$, let $\MT(n,r)$ be the set of all triangular arrays $M$ of the form
\setlength{\unitlength}{1.3mm}
\vs{-1}
\[\ba{c@{}c@{}c@{}c@{}c@{}c@{}c@{}c@{}c}
&&&M_{11}\,&\ldots&\;M_{1r}\\
&&M_{21}&&\ldots&&M_{2,2r}\\[-0.5mm]
&\bpic(1,1)\multiput(-1.5,-0.74)(1,1){3}{.}\epic&&&&&&\bpic(1,1)\multiput(-0.7,1.2)(1,-1){3}{.}\epic\\[1mm]
M_{n1}&&&&\ldots&&&&M_{n,nr}\ea\]
\vs{-8}

such that
\vs{-5}
\begin{itemize}
\item Each entry of $M$ is in $[n]$.\vs{-2}
\item In each row of $M$, any integer of $[n]$ appears at most $r$ times.\vs{-2}
\item $M_{ij}\le M_{i,j\pl1}$ \ for each $i\in[n]$, $j\in[ir\mi1]$.\vs{-2}
\item $M_{i\pl1,j}\le M_{ij}\le M_{i\pl1,j\pl r}$ \ for each $i\in[n\mi1]$, $j\in[ir]$.
\end{itemize}
\vs{-4}
It follows that the last row of any monotone triangle in $\MT(n,r)$ consists of each integer of $[n]$ repeated $r$ times.

It can be checked that there is a bijection between $\ASM(n,r)$ and $\MT(n,r)$
in which the monotone triangle $M$ which corresponds to the higher spin alternating sign matrix~$A$
is obtained by first using~(\ref{ASMToHV}) to find the vertical edge matrix~$V$ which corresponds to~$A$,
and then placing the integer $j$ $V_{ij}$ times in row $i$ of $M$, for each $i,j\in[n]$,
with these integers being placed in weakly increasing order along each row.
(Note that there is alternative bijection in which the horizontal edge matrix~$H$ which corresponds to~$A$
is obtained, and the integer $i$ is then placed $H_{ij}$ times in row $j$ of $M$, for each $i,j\in[n]$.)
For the inverse mapping, for each $i\in[0,n]$ and $j\in[n]$, $V_{ij}$ is set to be the number of times that $j$ occurs
in row $i$ of $M$, and $A$ is then obtained from $V$ using~(\ref{HVToASM}).

The set $\MT(n,1)$ of monotone triangles for standard alternating sign matrices
was introduced in~\cite{MilRobRum83}, and is also studied in, for example,~\cite{Fis06,Fis07,MilRobRum86,Pro01,Zei96a}.

The monotone triangle which corresponds to the running example of~(\ref{ExA}) and~(\ref{ExHV}) is
\begin{equation}\ba{cccccccccc}
&&&&2&3\\
&&&1&3&4&4\\
&&1&2&3&3&5&5\\
&1&1&2&3&3&4&5&5\\
1&1&2&2&3&3&4&4&5&5.\ea\end{equation}

Proceeding finally to sets of \emph{complementary edge matrix pairs}, for
$n\in\P$ and $r\in\N$ we define
\begin{equation}\label{CEM}\ba{l}\CEM(n,r)\,:=\\[2.5mm]
\;\;\left\{(\Hb,\Vb)\!=\!\left(\!
\left(\ba{c@{\;}c@{\;}c}\Hb_{10}&\ldots&\Hb_{1n}\\\vdots&&\vdots\\\Hb_{n0}&\ldots&\Hb_{nn}\ea\right)\!,\!
\left(\ba{c@{\;}c@{\;}c}\Vb_{01}&\ldots&\Vb_{0n}\\\vdots&&\vdots\\\Vb_{n1}&\ldots&\Vb_{nn}\ea\right)\!\right)\!
\in[0,r]^{n\times(n\pl1)}\times[0,r]^{(n\pl1)\times n}\:\right|\\[10.5mm]
\hs{8}\left.\ba{@{\bullet\ }l}
\Hb_{2k-1,0}=\Hb_{n-2k+2,n}=0, \ \ \Vb_{0,2k-1}=\Vb_{n,n-2k+2}=r,\mbox{ \ for all }k\in[\lceil\frac{n}{2}\rceil]\\[2mm]
\Hb_{2k,0}=\Hb_{n-2k+1,n}=r, \ \ \Vb_{0,2k}=\Vb_{n,n-2k+1}=0,\mbox{ \ for all }k\in[\lfloor\frac{n}{2}\rfloor]\\[2mm]
\Vb_{i\mi1,j}\pl\Hb_{i,j\mi1}\pl\Vb_{ij}\pl\Hb_{ij}=2r,\mbox{ \ for all }i,j\in[n]\ea\right\}\!.\!\!\!\ea\end{equation}

It can be seen that there is a bijection between $\EM(n,r)$ (and hence $\ASM(n,r)$) and $\CEM(n,r)$ in which
the complementary edge matrix pair $(\Hb,\Vb)$ which corresponds to the
edge matrix pair $(H,V)$ is given by
\begin{equation}\label{HVToHbVb}\ba{ll}\Hb_{ij}=\left\{\ba{l}
H_{ij}\,,\;i\pl j\mbox{ odd}\\[1.5mm]
r\mi H_{ij}\,,\;i\pl j\mbox{ even}\ea\right.\qquad&\mbox{for each }i\in[n],\ j\in[0,n]\\[8mm]
\Vb_{ij}=\left\{\ba{l}
r\mi V_{ij}\,,\;i\pl j\mbox{ odd}\\[1.5mm]
V_{ij}\,,\;i\pl j\mbox{ even}\ea\right.&\mbox{for each }i\in[0,n],\ j\in[n].\ea\end{equation}

The complementary edge matrix pair which corresponds to the running example of~(\ref{ExA}) and~(\ref{ExHV}) is
\begin{equation}\label{ExHbVb}(\Hb,\Vb)\,=\left(\left(\mbox{\footnotesize$\ba{cccccc}
0&2&1&0&2&0\\
2&1&2&0&0&2\\
0&2&1&0&0&0\\
2&1&1&1&0&2\\
0&2&1&1&2&0\ea$}\right),
\left(\mbox{\footnotesize$\ba{ccccc}
2&0&2&0&2\\
0&1&1&2&0\\
1&0&1&2&2\\
1&1&2&2&2\\
0&1&0&1&0\\
2&0&2&0&2\ea$}\right)\right).\end{equation}

\begin{figure}[h]
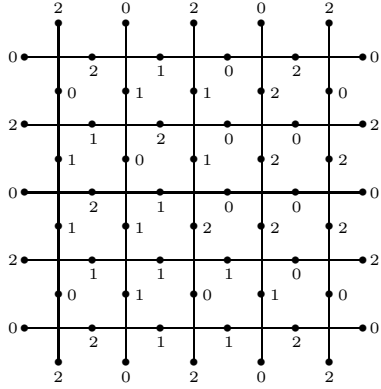
\centering
\bpic(6,6)\multiput(0.5,1)(0,1){5}{\line(1,0){5}}\multiput(1,0.5)(1,0){5}{\line(0,1){5}}
\multiput(0.5,1)(1,0){6}{\ppp{}{\bullet}}\multiput(0.5,2)(1,0){6}{\ppp{}{\bullet}}\multiput(0.5,3)(1,0){6}{\ppp{}{\bullet}}
\multiput(0.5,4)(1,0){6}{\ppp{}{\bullet}}\multiput(0.5,5)(1,0){6}{\ppp{}{\bullet}}
\multiput(1,0.5)(0,1){6}{\ppp{}{\bullet}}\multiput(2,0.5)(0,1){6}{\ppp{}{\bullet}}\multiput(3,0.5)(0,1){6}{\ppp{}{\bullet}}
\multiput(4,0.5)(0,1){6}{\ppp{}{\bullet}}\multiput(5,0.5)(0,1){6}{\ppp{}{\bullet}}
\multiput(1,5.65)(2,0){3}{\ppp{b}{2}}\multiput(1,0.35)(2,0){3}{\ppp{t}{2}}
\multiput(2,5.65)(2,0){2}{\ppp{b}{0}}\multiput(2,0.35)(2,0){2}{\ppp{t}{0}}
\multiput(0.4,1)(0,2){3}{\ppp{r}{0}}\multiput(5.6,1)(0,2){3}{\ppp{l}{0}}
\multiput(0.4,2)(0,2){2}{\ppp{r}{2}}\multiput(5.6,2)(0,2){2}{\ppp{l}{2}}
\put(1.13,4.48){\ppp{l}{0}}\put(2.13,4.48){\ppp{l}{1}}\put(3.13,4.48){\ppp{l}{1}}\put(4.13,4.48){\ppp{l}{2}}\put(5.13,4.48){\ppp{l}{0}}
\put(1.13,3.48){\ppp{l}{1}}\put(2.13,3.48){\ppp{l}{0}}\put(3.13,3.48){\ppp{l}{1}}\put(4.13,3.48){\ppp{l}{2}}\put(5.13,3.48){\ppp{l}{2}}
\put(1.13,2.48){\ppp{l}{1}}\put(2.13,2.48){\ppp{l}{1}}\put(3.13,2.48){\ppp{l}{2}}\put(4.13,2.48){\ppp{l}{2}}\put(5.13,2.48){\ppp{l}{2}}
\put(1.13,1.48){\ppp{l}{0}}\put(2.13,1.48){\ppp{l}{1}}\put(3.13,1.48){\ppp{l}{0}}\put(4.13,1.48){\ppp{l}{1}}\put(5.13,1.48){\ppp{l}{0}}
\put(1.52,4.87){\ppp{t}{2}}\put(1.52,3.87){\ppp{t}{1}}\put(1.52,2.87){\ppp{t}{2}}\put(1.52,1.87){\ppp{t}{1}}\put(1.52,0.87){\ppp{t}{2}}
\put(2.52,4.87){\ppp{t}{1}}\put(2.52,3.87){\ppp{t}{2}}\put(2.52,2.87){\ppp{t}{1}}\put(2.52,1.87){\ppp{t}{1}}\put(2.52,0.87){\ppp{t}{1}}
\put(3.52,4.87){\ppp{t}{0}}\put(3.52,3.87){\ppp{t}{0}}\put(3.52,2.87){\ppp{t}{0}}\put(3.52,1.87){\ppp{t}{1}}\put(3.52,0.87){\ppp{t}{1}}
\put(4.52,4.87){\ppp{t}{2}}\put(4.52,3.87){\ppp{t}{0}}\put(4.52,2.87){\ppp{t}{0}}\put(4.52,1.87){\ppp{t}{0}}\put(4.52,0.87){\ppp{t}{2}}
\epic\vs{-4}
\caption{\protect\ru{2}Assignment of entries of~(\ref{ExHbVb}) to lattice edges.\label{ExCEMC}}\end{figure}

In analogy with the association of an edge matrix pair to a configuration of a statistical mechanical model,
each entry of a complementary edge matrix pair can be assigned to an edge of the lattice,
i.e.,~$\Hb_{ij}$ is assigned to the horizontal edge between $(i,j)$ and $(i,j\pl1)$, for each $i\in[n]$, $j\in[0,n]$,
and~$\Vb_{ij}$ is assigned to the vertical edge between $(i,j)$ and $(i\pl1,j)$, for each $i\in[0,n]$, $j\in[n]$.
Also, in analogy with (\ref{Vr}), we define the set of \emph{complementary vertex types} as
\begin{equation}\VVb(r):=\{(\hb,\vb,\hb',\vb')\in[0,r]^4\mid \hb\pl\vb\pl\hb'\pl\vb'=2r\},\end{equation}
so that the lattice point~$(i,j)$
is associated with the complementary vertex type $(\Hb_{i,j\mi1},\Vb_{ij},\Hb_{ij},\Vb_{i\mi1,j})\in\VVb(r)$, for each $i,j\in[n]$.
Note that the mappings of each $(h,v,h',v')\in\V(r)$ to $(h,v,r\mi h',r\mi v')$,
or of each $(h,v,h',v')\in\V(r)$ to $(r\mi h,r\mi v,h',v')$, give two bijections between $\V(r)$ and $\VVb(r)$.
The assignment of the entries of the complementary edge matrix pair of~(\ref{ExHbVb}) to lattice edges is shown diagrammatically in
Figure~\ref{ExCEMC}.

\setlength{\unitlength}{9mm}

It is now natural to define, for each $n\in\P$ and $r\in\N$, the set $\FPL(n,r)$ of \emph{fully packed loop configurations}
as the set of all sets $P$ of nondirected open and closed lattice paths such that
\vs{-5}
\begin{itemize}
\item Successive points on each path of $P$ differ by $(-1,0)$, $(1,0)$, $(0,-1)$ or $(0,1)$.\vs{-2}
\item Each edge occupied by a path of $P$ is a horizontal edge between $(i,j)$ and $(i,j\pl1)$ with
$i\in[0,n]$ and $j\in[n]$, or a vertical edge between $(i,j)$ and $(i\pl1,j)$ with
$i\in[n]$ and $j\in[0,n]$.\vs{-2}
\item Any two edges occupied successively by a path of $P$ are different.\vs{-2}
\item Each edge is occupied by at most $r$ segments of paths of $P$.\vs{-2}
\item Each path of $P$ does not cross itself or any other path of $P$.\vs{-2}
\item Exactly $r$ segments of paths of $P$ pass through each (internal) point of $[n]\t[n]$.\vs{-2}
\item At each (external) point $(0,2k-1)$ and $(n\pl1,n-2k+2)$ for $k\in[\lceil\frac{n}{2}\rceil]$,
and $(2k,0)$ and $(n-2k+1,n\pl1)$ for $k\in[\lfloor\frac{n}{2}\rfloor]$, there are exactly $r$
endpoints of paths of $P$, these being the only lattice points which are path endpoints.
\end{itemize}
\vs{-4}
Note that an open nondirected lattice path is a sequence $(p_1,\ldots,p_m)$
of points of~$\Z^2$, for some $m\in\P$, where the reverse sequence $(p_m,\ldots,p_1)$ is regarded as the same path.
The endpoints of such a path are $p_1$ and $p_m$, and the pairs of successive points are $p_i$ and $p_{i\pl1}$, for
each $i\in[m\mi1]$.   A closed nondirected lattice path is a
sequence $(p_1,\ldots,p_m)$ of points of $\Z^2$, where
reversal and all cyclic permutations of the sequence are regarded as the same path.
Such a path has no endpoints, and its pairs of successive points are $p_i$ and $p_{i\pl1}$, for
each $i\in[m\mi1]$, as well as $p_1$ and~$p_m$.
For the case of $P\in\FPL(n,r)$, a path of $P$ whose points are all internal, i.e., in $[n]\t[n]$,
is closed, and a path of $P$ which has two external points, necessarily its endpoints, is open, even
if the two external points are the same.

\begin{figure}[h]
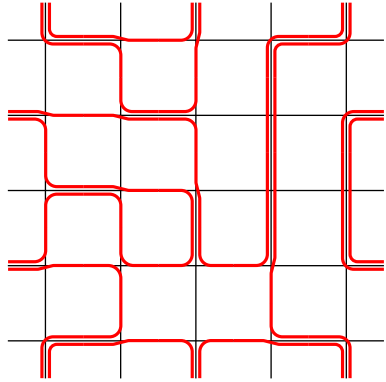
\centering
\pspicture(0.5,0.5)(5.5,5.7)
\multips(0.5,1)(0,1){5}{\psline[linewidth=0.5pt](0,0)(5,0)}\multips(1,0.5)(1,0){5}{\psline[linewidth=0.5pt](0,0)(0,5)}
\rput(1,5){\VVi}\rput(2,5){\VVxvii}\rput(3,5){\VVxi}\rput(4,5){\VVvi}\rput(5,5){\VVxiv}
\rput(1,4){\VVxvii}\rput(2,4){\VViv}\rput(3,4){\VVxviii}\rput(4,4){\VVviii}\rput(5,4){\VVvi}
\rput(1,3){\VVii}\rput(2,3){\VVxvii}\rput(3,3){\VVxvi}\rput(4,3){\VVviii}\rput(5,3){\VVviii}
\rput(1,2){\VVxiii}\rput(2,2){\VVxb}\rput(3,2){\VVv}\rput(4,2){\VVxi}\rput(5,2){\VVi}
\rput(1,1){\VVvi}\rput(2,1){\VVxiii}\rput(3,1){\VVxv}\rput(4,1){\VViv}\rput(5,1){\VVxix}
\endpspicture
\caption{\protect\ru{2}A fully packed loop configuration which maps to~(\ref{ExHbVb}).\label{ExFPL}}\end{figure}

It can now be seen that there is a mapping from $\FPL(n,r)$ to $\CEM(n,r)$
in which the fully packed loop configuration $P$ is mapped to the complementary edge matrix pair $(\Hb,\Vb)$ according to
\begin{equation}\label{PToHbVb}\ba{r@{\;\;}l}\Hb_{ij}\;=&\mbox{number of segments of paths of }P\mbox{ which occupy the edge between}\\[1mm]
&(i,j)\mbox{ and }(i,j\pl1),\mbox{ \ for each }i\in[n],\,j\in[0,n]\\[3mm]
\Vb_{ij}\;=&\mbox{number of segments of paths of }P\mbox{ which occupy the edge between }\\[1mm]
&(i\pl1,j)\mbox{ and }(i,j),\mbox{ \ for each }i\in[0,n],\,j\in[n].\ea\end{equation}
A fully packed loop configuration of $\FPL(5,2)$ which maps to the complementary edge matrix pair of~(\ref{ExHbVb})
is shown diagrammatically in Figure~\ref{ExFPL}.

\psset{unit=10mm}

It can be checked that the mapping of~(\ref{PToHbVb}) is surjective for each $r\in\N$ and $n\in\P$.
Furthermore, for $r\in\{0,1\}$ or $n\in\{1,2\}$ it is injective, while for $r\ge2$ and $n\ge3$ it is not injective.
This is due to the fact that if, for a complementary vertex type $(\Hb_{i,j\mi1},\Vb_{ij},\Hb_{ij},\Vb_{i\mi1,j})\in\VVb(r)$,~(\ref{PToHbVb})
is used to assign appropriate numbers of path segments to the four edges surrounding the
point $(i,j)\in[n]\t[n]$, then for $r\in\{0,1\}$ there is always a unique way to link these segments through $(i,j)$,
whereas for $r\ge2$ there can be several ways of linking these segments through $(i,j)$, such cases occurring for each $n\ge3$.
For example, for $r=2$ there is a unique way of linking the segments, except if
$(\Hb_{i,j\mi1},\Vb_{ij},\Vb_{i\mi1,j},\Hb_{ij})=(1,1,1,1)$, in which case either of the
configurations
\raisebox{-2.6ex}{\pspicture(0,-0.1)(1,1.1)\Vi\rput(0.5,0.5){\VVxa}\endpspicture}
or
\raisebox{-2.6ex}{\pspicture(0,-0.1)(1,1.1)\Vi\rput(0.5,0.5){\VVxb}\endpspicture} can be used.
Thus, since the example $(\Hb,\Vb)$ of~(\ref{ExHbVb}) and Figure~\ref{ExCEMC} has the single case $(i,j)=(2,2)$ where this occurs, there
are two fully packed loop configurations of $\FPL(5,2)$ which map to $(\Hb,\Vb)$: that of Figure~\ref{ExFPL}
and that which differs from it by the configuration at $(2,2)$.

The cases of $\FPL(n,1)$, and of certain related sets which arise by imposing additional symmetry conditions,
have been studied extensively.  See for example~\cite{CasKra04,CasKraLasNad05,DifZinZub04,DifZub04,Wie00,Zub04}.
In these studies, the fully packed loop configurations~$P$ of $\FPL(n,1)$ are usually classified according to the
\emph{link pattern} formed among the external points by the open paths of~$P$.  This then leads to important
results and conjectures, including unexpected connections with certain statistical mechanical models.
See for example~\cite{Deg05,Deg07} and references therein.

Link patterns for certain higher spin statistical mechanical models have been studied in~\cite{Zin07}, although
an additional condition is imposed there, which for fully packed loop configurations would be that
the two endpoints of any open path of $P\in\FPL(n,r)$ must be distinct.  However, the subset of
$\FPL(n,r)$ in which this condition is satisfied is still not in bijection with $\CEM(n,r)$, for $n\ge3$ and $r\ge2$ .
For example, the fully packed loop configurations
\psset{unit=8mm}
\raisebox{-6.6ex}{
\pspicture(0.5,0.3)(3.7,3.6)
\multips(0.5,1)(0,1){3}{\psline[linewidth=0.5pt](0,0)(3,0)}\multips(1,0.5)(1,0){3}{\psline[linewidth=0.5pt](0,0)(0,3)}
\rput(1,3){\VViii}\rput(2,3){\VVix}\rput(3,3){\VVxiv}
\rput(1,2){\VVxiii}\rput(2,2){\VVxb}\rput(3,2){\VVix}
\rput(1,1){\VVvi}\rput(2,1){\VVxiii}\rput(3,1){\VVxvi}
\endpspicture} and
\raisebox{-6.6ex}{
\pspicture(0.5,0.3)(3.7,3.6)
\multips(0.5,1)(0,1){3}{\psline[linewidth=0.5pt](0,0)(3,0)}\multips(1,0.5)(1,0){3}{\psline[linewidth=0.5pt](0,0)(0,3)}
\rput(1,3){\VViii}\rput(2,3){\VVix}\rput(3,3){\VVxiv}
\rput(1,2){\VVxiii}\rput(2,2){\VVxa}\rput(3,2){\VVix}
\rput(1,1){\VVvi}\rput(2,1){\VVxiii}\rput(3,1){\VVxvi}
\endpspicture}
in $\FPL(3,2)$
both satisfy the condition and map to the same complementary edge matrix pair in $\CEM(3,2)$ using~(\ref{PToHbVb}).

\subsubsection{The Alternating Sign Matrix Polytope}
In this section, we define the alternating sign matrix polytope in $\R^{n^2}$, using a halfspace description, and we show that
its vertices are the standard alternating sign matrices of size $n$.

We begin by summarizing the facts about convex polytopes which will be needed here.
For further information, see for example~\cite{Zie95}.
For $m\in\P$, a convex polytope in $\R^m$ can be defined as
a bounded intersection of finitely-many closed affine halfspaces in $\R^m$,
or equivalently as a convex hull of finitely-many points in~$\R^m$.  The equivalence
of these descriptions is nontrivial and is proved, for example, in~\cite[Theorem~1.1]{Zie95}.
It follows that hyperplanes in $\R^m$ can be included together with closed
halfspaces in the first description, since a hyperplane is simply the intersection of the two closed halfspaces which
meet at the hyperplane.
The dimension, $\dim\PP$, of a convex polytope $\PP$ is defined to be the dimension of its affine hull,
$\aff(\PP):=\{\l p_1\pl(1\mi\l)p_2\mid p_1,p_2\in\PP,\ \l\in\R\}$.
A face of $\PP$ is an intersection of $\PP$ with any hyperplane for which $\PP$ is a subset of one of the two
closed halfspaces determined by the hyperplane.  If a face contains only one point, that point is known as a vertex.
Thus, the set of vertices, $\ver\PP$, of $\PP\subset\R^m$ is the set of points $p\in\PP$
for which there exists a closed affine halfspace $\S$ in $\R^m$ such that $\PP\cap\S=\{p\}$.
It can be shown that $\ver\PP$ is also the set of points $p\in\PP$ which do not lie in the interior of any line segment in $\PP$,
i.e., $p\in\PP$ is a vertex of $\PP$ if and only if there do not exist $\l\in\R_{(0,1)}$ and $p_1\ne p_2\in\PP$ with
$p=\l p_1+(1\mi\l)p_2$.  Any convex polytope~$\PP$ has only finitely-many vertices,
and is the convex hull of these vertices, or equivalently the set of all convex combinations of its vertices,
$\PP=\{\sum_{p\in\ver\PP}\l_p\:p\mid\l_p\in\R_{[0,1]}\mbox{ for each }p\in\ver\PP,\ \sum_{p\in\ver\PP}\l_p=1\}$.
Convex polytopes $\PP\subset\R^m$ and $\PP'\subset\R^{m'}$
are defined to be
affinely isomorphic if there is an affine map $\phi:\R^m\rightarrow\R^{m'}$ which is bijective between
$\PP$ and $\PP'$. In such cases,
$\ver\PP'=\phi(\ver\PP)$.  Finally, a convex polytope whose vertices all have integer coordinates is known as integral.

We now define, for $n\in\P$,
\begin{equation}\label{An}\ba{l}\!\!\!\A_n:=\\[2mm]
\hs{4.4}\left\{\!x\!=\!\!\left(\ba{c@{\;}c@{\;}c}
x_{11}&\ldots&x_{1n}\\\vdots&&\vdots\\x_{n1}&\ldots&x_{nn}\ea\right)\!\!\in\R^{n\times n}\:\left|\;\,\ba{l}
\bullet\ \sum_{j'=1}^n\!x_{ij'}=\sum_{i'=1}^n\!x_{i'\!j}=1\mbox{ \ for all }i,j\in[n]\!\\[2.2mm]
\bullet\ \sum_{j'=1}^jx_{ij'\!}\ge0\mbox{ \ for all }i\in[n],\;j\in[n\mi1]\\[2.2mm]
\bullet\ \sum_{j'=j}^nx_{ij'\!}\ge0\mbox{ \ for all }i\in[n],\;j\in[2,n]\\[2.2mm]
\bullet\ \sum_{i'=1}^ix_{i'\!j}\ge0\mbox{ \ for all }i\in[n\mi1],\;j\in[n]\\[2.2mm]
\bullet\ \sum_{i'=i}^nx_{i'\!j}\ge0\mbox{ \ for all }i\in[2,n],\;j\in[n]\ea\right.\right\}\\[20mm]
=\left\{\!x\!=\!\!\left(\ba{c@{\;}c@{\;}c}
x_{11}&\ldots&x_{1n}\\\vdots&&\vdots\\x_{n1}&\ldots&x_{nn}\ea\right)\!\!\in\R^{n\times n}\:\left|\;\ba{l}
\bullet\ \sum_{j'=1}^nx_{ij'}=\sum_{i'=1}^nx_{i'\!j}=1\mbox{ \ for all }i,j\in[n]\\[2.2mm]
\bullet\ 0\le\sum_{j'=1}^jx_{ij'\!}\le1\mbox{ \ for all }i\in[n],\;j\in[n\mi1]\\[2.2mm]
\bullet\ 0\le\sum_{i'=1}^ix_{i'\!j}\le1\mbox{ \ for all }i\in[n\mi1],\;j\in[n]\ea\right.\!\right\}\!.\ea\!\!\!\!\!\!\!\end{equation}
In other words, $\A_n$ is the set of $n\t n$ real-entry matrices
for which all complete row and column sums are $1$,
and all partial row and column sums extending from each end of the row or column are nonnegative.

It can be seen that each entry of any matrix of $\A_n$ is between $-1$ and~$1$,
and that if the entry is in the first or last row or column, then it is
between $0$ and $1$, so that~$\A_n$ is a bounded subset of~$\R^{n^2}$.
Since $\A_n$ is also an intersection of finitely-many closed halfspaces and
hyperplanes in~$\R^{n^2}$, it is a convex polytope in $\R^{n^2}$, and will be referred to as
the \emph{alternating sign matrix polytope}. This polytope was defined independently,
using a convex hull description, in~\cite{Str07}.

An example of an element of $\A_4$ is
\begin{equation}\label{Exx}x\:=\left(\mbox{\footnotesize$
\ba{ccccc}.3&0&.6&.1\\
.2&.5&-.6&.9\\
.5&-.5&1&0\\
0&1&0&0\ea$}\right).\end{equation}

Defining
\begin{equation}\B_n:=\{x\in\A_n\mid x_{ij}\ge0\mbox{ \ for each }i,j\in[n]\}\,,\end{equation}
it can be seen that this is the set of \emph{doubly stochastic matrices} of size $n$, i.e., nonnegative real-entry
$n\times n$ matrices for which all complete row and column sums are~$1$.  This is the convex polytope in $\R^{n^2}$
often known as the \emph{Birkhoff polytope} (or \emph{assignment polytope}).
See for example~\cite[Ex.~0.12]{Zie95} and references therein.

It now follows that the higher spin alternating sign matrices and semimagic squares of size $n$
with line sum $r$ are the integer points of the $r$-th dilates of $\A_n$ and $\B_n$ respectively,
\begin{equation}\label{dil}
\ASM(n,r)\,=\,r\A_n\cap\Z^{n^2}\,,\qquad
\SMS(n,r)\,=\,r\B_n\cap\N^{n^2},
\end{equation}
where the $r$-th dilate of a set $P\subset\R^m$ is simply $r P:=\{r x\mid x\in P\}$.

It also follows that $\aff\A_n=\aff\B_n=
\{x\in\R^{n\times n}\mid\sum_{j'=1}^n\!x_{ij'}=\sum_{i'=1}^n\!x_{i'\!j}=1\mbox{ for all }$ $i,j\in[n]\}$,
and that of the $2n$ linear equations in $n^2$ variables within this set, only
$2n\mi1$ equations are independent, so that
\begin{equation}\label{dimA}\dim\A_n=\dim\B_n=(n\mi1)^2.\end{equation}
This is effectively equivalent to the fact that any $x\in\aff\A_n=\aff\B_n$
can be obtained by freely choosing the entries of any $(n\mi1)\t(n\mi1)$ submatrix of $x$,
the remaining $2n\mi1$ entries of $x$ then being determined by the condition that each row and column
sum is~$1$.

We also define, for $n\in\P$,
\begin{equation}\label{E}\ba{l}\E_n:=
\Biggl\{(h,v)\!=\!\left(\!
\left(\ba{c@{\;}c@{\;}c}h_{10}&\ldots&h_{1n}\\\vdots&&\vdots\\h_{n0}&\ldots&h_{nn}\ea\right)\!,\!
\left(\ba{c@{\;}c@{\;}c}v_{01}&\ldots&v_{0n}\\\vdots&&\vdots\\v_{n1}&\ldots&v_{nn}\ea\right)\!\right)\!
\in\R_{[0,1]}{}^{n\times(n\pl1)}\times\R_{[0,1]}\!^{(n\pl1)\times n}\:\Bigg|\\[8mm]
\hs{12}h_{i0}=v_{0j}=0,\ h_{in}=v_{nj}=1,\ h_{i,j\mi1}\pl v_{ij}=v_{i\mi1,j}\pl h_{ij},
\mbox{ \ for all }i,j\in[n]\!\Biggr\}.\!\!\!\ea\end{equation}
This is a convex polytope in $\R^{2n(n\pl1)}$, which we shall refer to as the \emph{edge matrix polytope}.
It can be seen that $\EM(n,r)\,=\,r\E_n\cap\Z^{2n(n\pl1)}$, and that, analogously to~(\ref{ASMToHV}) and~(\ref{HVToASM}),
there is a bijection between $\A_n$ and $\E_n$ in which the $(h,v)\in\E_n$ which corresponds to $x\in\A_n$ is given by
\begin{equation}\label{AToE}\ba{c}\ds h_{ij}=\sum_{j'=1}^jx_{ij'}\,,\mbox{ \ for each }i\in[n],\ j\in[0,n]\\[5mm]
\ds v_{ij}=\sum_{i'=1}^ix_{i'j}\,,\mbox{ \ for each }i\in[0,n],\ j\in[n],\ea\end{equation}
and inversely,
\begin{equation}\label{EToA}x_{ij}\,=\,h_{ij}-h_{i,j\mi1}\,=\,v_{ij}-v_{i\mi1,j}\,,\mbox{ \ for each }i,j\in[n].\end{equation}
Furthermore,~(\ref{AToE}) can be extended to a linear map from $\R^{n^2}$ to $\R^{2n(n\pl1)}$, or~(\ref{EToA}) can be extended to
a linear map from $\R^{2n(n\pl1)}$ to $\R^{n^2}$, implying that $\A_n$ and $\E_n$ are affinely isomorphic.
Also, in analogy with~(\ref{HVSum}),
\begin{equation}\label{hvSum}\sum_{i,j=1}^nh_{ij}\,=\,\sum_{i,j=1}^nv_{ij}\,=\,n(n\pl1)/2\mbox{ \ for each }(h,v)\in\E_n.\end{equation}

It is shown in~\cite{Bir46,Von53} that the vertices of the Birkhoff polytope $\B_n$ are the
permutation matrices of size~$n$, so that $\B_n$ is an integral convex polytope.
We now state and prove the corresponding result for~$\A_n$.
This result was obtained independently in~\cite{Str07}.
\begin{theorem}
The vertices of the alternating sign matrix polytope $\A_n$ are the
standard alternating sign matrices of size $n$.
\end{theorem}
\textit{Proof.} \ We shall show that the vertices of the edge matrix polytope $\E_n$
are the edge matrix pairs of $\EM(n,1)$, i.e., $\ver\E_n=\EM(n,1)$.  It then follows, since~$\E_n$ and $\A_n$ are
affinely isomorphic with mapping~(\ref{EToA}),
and since~(\ref{EToA}) maps $\EM(n,1)$ to $\ASM(n,1)$, that $\ver\A_n=\ASM(n,1)$ as required.

We first show that $\EM(n,1)\subset\ver\E_n$.  Consider any $(H,V)\in\EM(n,1)$.
From~(\ref{EM}), this is a pair of matrices with a total of $2n(n\pl1)$ $\{0,1\}$-entries, of which, due to~(\ref{HVSum}),
$n(n\pl1)$ are $0$'s and $n(n\pl1)$ are $1$'s.
Now define the halfspace $\S=\{(y,z)\in\R^{n\times(n\pl1)}\t\R^{(n\pl1)\times n}\mid
\sum_{i,j=1}^{n}(H_{ij}\,y_{ij}+V_{ij}\,z_{ij})\ge n(n\pl1)\}$,
and consider any matrix pair $(h,v)\in\E_n\cap\S$.
Due to~(\ref{EM}),~(\ref{E}) and~(\ref{hvSum}), one such matrix pair is $(H,V)$.
Also, $(h,v)\in\E_n$ implies, using~(\ref{E}) and~(\ref{hvSum}), that
each of the $2n(n\pl1)$ entries of $(h,v)$ is between $0$ and $1$ inclusive, and that they all sum to $n(n\pl1)$,
while $(h,v)\in\S$ implies that the $n(n\pl1)$ entries of $(h,v)$ in
the same positions as the~$1$'s of $(H,V)$ sum to at least $n(n\pl1)$.
It can be seen that these conditions are only satisfied if $(h,v)=(H,V)$.
Therefore, $\E_n\cap\S=\{(H,V)\}$, implying that $(H,V)\in\ver\E_n$ as required.
(Note that alternatively it could have been shown here that
$\E_n\cap\{(y,z)\in\R^{n\times(n\pl1)}\t\R^{(n\pl1)\times n}\mid\sum_{i,j=1}^{n}H_{ij}\,y_{ij}\ge n(n\pl1)/2\}=\{(H,V)\}$
or that $\E_n\cap\{(y,z)\in\R^{n\times(n\pl1)}\t\R^{(n\pl1)\times n}\mid\sum_{i,j=1}^{n}V_{ij}\,z_{ij}\ge n(n\pl1)/2\}=\{(H,V)\}$.)

We now show that $\ver\E_n\subset\EM(n,1)$.
Consider any $(h,v)\in\E_n\!\setminus\EM(n,1)$.
We shall eventually deduce that $(h,v)\notin\ver\E_n$, which gives the required result.
Similarly to the association of edge matrix pairs with configurations of a statistical mechanical model,
we associate~$h_{ij}$ with the horizontal edge between lattice points $(i,j)$ and $(i,j\pl1)$, for each $i\in[n]$, $j\in[0,n]$,
and~$v_{ij}$ with the vertical edge between lattice points $(i,j)$ and $(i\pl1,j)$, for each $i\in[0,n]$, $j\in[n]$
(using matrix-type labeling of lattice points).  Since $\EM(n,1)=\E_n\cap\Z^{2n(n\pl1)}$, $(h,v)\notin\EM(n,1)$ implies
that at least one entry of $(h,v)$ is nonintegral.  Now, from~(\ref{E}), $(h,v)\in\E_n$ implies
that $h_{i0}=v_{0j}=0$ and $h_{in}=v_{nj}=1$, for each $i,j\in[n]$,
so that any nonintegral entry of $(h,v)$ must be associated with one of the $2n(n\mi1)$ internal edges,
(i.e., the horizontal edges between $(i,j)$ and $(i,j\pl1)$, for $i\in[n]$, $j\in[n\mi1]$, and the vertical
edges between $(i,j)$ and $(i\pl1,j)$, for each $i\in[n\mi1]$, $j\in[n]$).
Also from~(\ref{E}),~$(h,v)\in\E_n$ implies,  that
\begin{equation}\label{hvCond}h_{i,j\mi1}\pl v_{ij}\,=\,v_{i\mi1,j}\pl h_{ij}\,,\end{equation}
for each $i,j\in[n]$.  But if any one of the four
entries in~(\ref{hvCond}) is nonintegral, then at least one of the others must also be nonintegral.
Therefore, the existence among the entries of
$(h,v)$ of a noninteger implies the existence of two further nonintegers, among each of the other three entries
of the two cases of~(\ref{hvCond}) in which the initial nonintegral entry appears.
It now follows by repeatedly applying this argument,
and since the internal edges form a finite and closed grid,
that there exists at least one cycle of internal
edges associated with noninteger entries of $(h,v)$.

\begin{figure}[h]
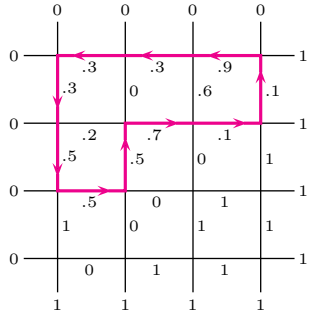
\centering
\pspicture(0.5,0)(4.5,4.7)
\multips(0.5,1)(0,1){4}{\psline[linewidth=0.5pt](0,0)(4,0)}\multips(1,0.5)(1,0){4}{\psline[linewidth=0.5pt](0,0)(0,4)}
\multirput[b](1,4.6)(1,0){4}{$\sss0$}\multirput[t](1,0.4)(1,0){4}{$\sss1$}\multirput[r](0.43,1)(0,1){4}{$\sss0$}\multirput[l](4.56,1)(0,1){4}{$\sss1$}
\rput[l](1.06,3.52){$\sss.3$}\rput[l](2.06,3.48){$\sss0$}\rput[l](3.06,3.48){$\sss.6$}\rput[l](4.06,3.48){$\sss.1$}
\rput[l](1.06,2.52){$\sss.5$}\rput[l](2.06,2.48){$\sss.5$}\rput[l](3.06,2.48){$\sss0$}\rput[l](4.06,2.48){$\sss1$}
\rput[l](1.06,1.48){$\sss1$}\rput[l](2.06,1.48){$\sss0$}\rput[l](3.06,1.48){$\sss1$}\rput[l](4.06,1.48){$\sss1$}
\rput[t](1.47,3.91){$\sss.3$}\rput[t](1.47,2.91){$\sss.2$}\rput[t](1.47,1.91){$\sss.5$}\rput[t](1.47,0.91){$\sss0$}
\rput[t](2.47,3.91){$\sss.3$}\rput[t](2.43,2.91){$\sss.7$}\rput[t](2.47,1.91){$\sss0$}\rput[t](2.47,0.91){$\sss1$}
\rput[t](3.47,3.91){$\sss.9$}\rput[t](3.47,2.91){$\sss.1$}\rput[t](3.47,1.91){$\sss1$}\rput[t](3.47,0.91){$\sss1$}
\psline[linewidth=1.2pt,linecolor=magenta]{->}(4,3.75)(4,4)(3.2,4)
\psline[linewidth=1.2pt,linecolor=magenta]{->}(3.25,4)(2.2,4)
\psline[linewidth=1.2pt,linecolor=magenta]{->}(2.25,4)(1.2,4)
\psline[linewidth=1.2pt,linecolor=magenta]{->}(1.25,4)(1,4)(1,3.2)
\psline[linewidth=1.2pt,linecolor=magenta]{->}(1,3.25)(1,2.2)
\psline[linewidth=1.2pt,linecolor=magenta]{->}(1,2.25)(1,2)(1.8,2)
\psline[linewidth=1.2pt,linecolor=magenta]{->}(1.75,2)(2,2)(2,2.8)
\psline[linewidth=1.2pt,linecolor=magenta]{->}(2,2.75)(2,3)(2.8,3)
\psline[linewidth=1.2pt,linecolor=magenta]{->}(2.75,3)(3.8,3)
\psline[linewidth=1.2pt,linecolor=magenta]{->}(3.75,3)(4,3)(4,3.8)
\endpspicture\vs{-4}
\caption{A cycle of nonintegers for the example of (\ref{Exx}).\label{ExNI}}\end{figure}

We select any such cycle, give it an orientation, say anticlockwise, and denote the sets of
points $(i,j)$ for which the horizontal edge between $(i,j)$ and $(i,j\pl1)$ is in the cycle and directed right
or left as respectively~$\HH_+$ or~$\HH_-$, and the sets of
points $(i,j)$ for which the vertical edge between $(i,j)$ and $(i\pl1,j)$ is in the cycle
and directed up or down as respectively~$\V_+$ or~$\V_{\!-}$.  An example of such a cycle,
for the $(h,v)\in\E_4$ which corresponds to the example $x\in\A_4$ of~(\ref{Exx}) is shown diagrammatically in Figure~\ref{ExNI}.
\psset{unit=9mm}
For this example, $\HH_+=\{(2,2),\,(2,3),\,(3,1)\}$, $\HH_-=\{(1,1),\,(1,2),\,(1,3)\}$,
$\V_{\!+}=\{(1,4),\,(2,2)\}$ and $\V_{\!-}=\{(1,1),\,(2,1)\}$.
We now define, for $\m\in\R$, a matrix pair $(h'(\m),v'(\m))$, with entries $h'(\m)_{ij}$ for $i\in[n]$, $j\in[0,n]$, and $v'(\m)_{ij}$
for $i\in[0,n]$, $j\in[n]$, given by
\begin{equation}\label{hvm}h'(\m)_{ij}=\left\{\ba{l}
h_{ij}\pl\m\,,\;(i,j)\in\HH_+\\[1.5mm]
h_{ij}\mi\m\,,\;(i,j)\in\HH_-\\[1.5mm]
h_{ij}\,,\;\mathrm{otherwise}\ea\right.\qquad\qquad
v'(\m)_{ij}=\left\{\ba{l}
v_{ij}\pl\m\,,\;(i,j)\in\V_{\!+}\\[1.5mm]
v_{ij}\mi\m\,,\;(i,j)\in\V_{\!-}\\[1.5mm]
v_{ij}\,,\;\mathrm{otherwise}\,.\ea\right.\end{equation}
For the example of Figure~\ref{ExNI},
\begin{equation}(h(\m),v(\m))\,=\left(\left(\mbox{\footnotesize$\ba{ccccc}
0&.3\mi\m&.3\mi\m&.9\mi\m&1\\
0&.2&.7\pl\m&.1\pl\m&1\\
0&.5\pl\m&0&1&1\\
0&0&1&1&1\ea$}\right),\;\left(\mbox{\footnotesize$\ba{cccc}0&0&0&0\\
.3\mi\m&0&.6&.1\pl\m\\
.5\mi\m&.5\pl\m&0&1\\
1&0&1&1\\
1&1&1&1\ea$}\right)\right)\!.\end{equation}
We now check whether $(h'(\m),v'(\m))\in\E_n$.
By using~(\ref{hvm}) to replace each entry of $(h,v)$ in~(\ref{hvCond}) with an entry of
$(h'(\m),v'(\m))$, it can be checked that
the equation $h'(\m)_{i,j\mi1}\pl v'(\m)_{ij}=v'(\m)_{i\mi1,j}\pl h'(\m)_{ij}$ is
satisfied for each $i,j\in[n]$, since if the cycle does not pass through $(i,j)$,
then the required equation is immediately obtained, while
for all possible configurations (of which there are six),
and both possible directions, in which the cycle can pass through~$(i,j)$, all explicit appearances of~$\m$ cancel out,
again leaving the required equation.
The conditions
$h'(\m)_{i0}=v'(\m)_{0j}=0$ and $h'(\m)_{in}=v'(\m)_{nj}=1$
are also met for each $i,j\in[n]$, since $(h'(\m),v'(\m))$ and $(h,v)$ match for these entries.
It only remains for the conditions $h'(\m)_{ij},v'(\m)_{ij}\in\R_{[0,1]}$
to be checked for each $i,j\in[n]$, and it can be seen that these
are satisfied if $-\m_-\le\m\le\m_+$, where
\begin{equation}\ba{r@{\;}c@{\;}l}\m_-&:=&\min(\{h_{ij}\mid(i,j)\in\HH_+\}\cup\{v_{ij}\mid(i,j)\in\V_{\!+}\}\\[1.5mm]
&&\qquad\quad\mbox{}\cup\{1\mi h_{ij}\mid(i,j)\in\HH_-\}\cup\{1\mi v_{ij}\mid(i,j)\in\V_{\!-}\})\\[3mm]
\m_+&:=&\min(\{1\mi h_{ij}\mid(i,j)\in\HH_+\}\cup\{1\mi v_{ij}\mid(i,j)\in\V_{\!+}\}\\[1.5mm]
&&\qquad\quad\mbox{}\cup\{h_{ij}\mid(i,j)\in\HH_-\}\cup\{v_{ij}\mid(i,j)\in\V_{\!-}\})\,.\ea\end{equation}
The facts that $h_{ij}\in\R_{(0,1)}$ for each $(i,j)\in\HH_\pm$, and $v_{ij}\in\R_{(0,1)}$ for each $(i,j)\in\V_{\!\pm}$,
imply that $\m_-$ and $\m_+$ are both positive, so that there exists a finite-length, closed interval of values of $\m$ for
which $(h'(\m),v'(\m))\in\E_n$.
For the example of Figure~\ref{ExNI}, $\m_-=.1$ and $\m_+=.3$.
Finally, we define
$(h_{\pm},v_{\pm}):=(h'(\pm\m_\pm),v'(\pm\m_\pm))$,
so that it follows that $(h_\pm,v_\pm)\in\E_n$, $(h_-,v_-)\ne(h_+,v_+)$, and
\begin{equation}\label{LI}\ts(h,v)=\frac{\m_+}{\m_-+\m_+}(h_-,v_-)+\frac{\m_-}{\m_-+\m_+}(h_+,v_+)\,.\end{equation}
Therefore, $(h,v)$ lies in the interior of a line segment between two points of $\E_n$, and so, as required, $(h,v)\notin\ver\E_n$.
This concludes the proof of Theorem~1.\hspace{\fill}$\Box$

It follows immediately from Theorem~1 that the alternating sign matrix polytope is integral, a fact which
will be important in the enumeration of higher spin alternating sign matrices in the next section.

It also follows from Theorem~1 that $\A_n$ can be
described as the convex hull of the alternating sign matrices of size $n$, or equivalently that
\begin{equation}\label{AD}\ba{l}\A_n=\\[1.5mm]
\hs{7}\ds\Biggl\{\sum_{A\in\ASM(n,1)}\!\!\l_A\:A\;\Bigg|\;
\l_A\in\R_{[0,1]}\mbox{ \ for each }A\in\ASM(n,1),\ \sum_{A\in\ASM(n,1)}\!\!\l_A=1\Biggr\}.\ea\!\!
\end{equation}
The way in which this conclusion has been reached here depends on the general theorem that any
convex polytope has both a halfspace and convex hull description.  However,~(\ref{AD}) can also be derived more directly.
It follows immediately from~(\ref{ASM}) and~(\ref{An}) that the RHS of~(\ref{AD}) is a subset of the LHS, i.e.,
that every convex combination of standard alternating sign matrices of size $n$ is an element of~$\A_n$.
Conversely, every element of $\A_n$ is a convex combination of elements of $\ASM(n,1)$
due to the following argument, which we owe to~\cite{Str07}.  For each $x\in\A_n$, let $\nu(x)$ be the number of nonintegers
among the entries of $(h,v)$, where $(h,v)\in\E_n$ corresponds to~$x$ by~(\ref{AToE}). Now consider
the following recursive process applied to $x\in\A_n$.
If $\nu(x)=0$, then $x\in\ASM(n,1)$, and the process is terminated.
If $\nu(x)$ is positive,
then let $x_\pm\in\A_n$ correspond to $(h_\pm,v_\pm)\in\E_n$,
where $(h_\pm,v_\pm)$ are defined in terms of $(h,v)$ as in the proof of Theorem~1.
Using~(\ref{LI}), write $x$ as the convex combination
$x=\frac{\m_+}{\m_-+\m_+}x_-+\frac{\m_-}{\m_-+\m_+}x_+$, where it can be seen
that $\nu(x_\pm)<\nu(x)$.  Now reapply the process to $x_-$ and $x_+$.
It follows that any $x\in\A_n$ can thus be
completely decomposed as a convex combination of elements of $\ASM(n,1)$.

Faces of the alternating sign polytope other than
those given by its vertices are studied in detail in~\cite{Str07}, but we shall not consider these here.

\subsubsection{Enumeration of Higher Spin Alternating Sign Matrices of Fixed Size}
In this section, we use the general theory of the enumeration of integer points in
integer dilates of integral convex polytopes to obtain results on the
enumeration of higher spin alternating sign matrices of fixed size.

We begin by summarizing the details of this theory which will be needed here.
For further information, more general results, and references, see for example~\cite{BecRob07} or~\cite[Sec.~4.6]{Sta86}.
For an integral convex polytope $\PP$ in $\R^m$ with relative interior $\PO$,
there exists a unique function $\Lp(r)$,
where $r$ can be regarded as a complex variable, with the properties that:
\begin{equation}\label{Ehr}
\ba{r@{\;\;}l}
\mathrm{(i)}&\Lp(r)\mbox{ is a polynomial in $r$ of degree }\dim\PP.\\[3mm]
\mathrm{(ii)}&|r\PP\cap\Z^m|\,=\,\Lp(r),\mbox{ \ for each }r\in\N.\\[3mm]
\mathrm{(iii)}&|r\PO\cap\Z^m|\,=\,(-1)^{\dim\PP}\,\Lp(-r),\mbox{ \ for each }r\in\N.\\[3mm]
\mathrm{(iv)}&\Lp(r)\mbox{ can be expressed in the form \ }
\ds\Lp(r)\,=\sum_{k=0}^{\dim\PP}c_k\left(\,\ba{c}r+k\\[1.2mm]\dim\PP\ea\,\right),\\[5mm]
&\mbox{with }c_k\in\N\mbox{ for each }k\in[0,\dim\PP].
\ea\end{equation}
The function $\Lp(r)$ is known as the \emph{Ehrhart polynomial} of $\PP$.  It
can be seen that $\Lp(0)=c_{\dim\PP}=1$, and that
the Ehrhart polynomial can be obtained explicitly by finding $\Lp(r)$ for $\dim\PP$ further integers $r$
by directly enumerating the number of lattice points $|r\PP\cap\Z^m|$ for positive $r$, or $|-r\PO\cap\Z^m|$
for negative $r$, and then interpolating.

For the alternating sign matrix polytope $\A_n$, the relative interior $\A_n^\circ$ is obtained
by simply replacing each weak inequality in~(\ref{An}) by a strict inequality.  Defining, for $n\in\P$
and $r\in\N$,
\begin{equation}\label{ASMOdil}\ASMO(n,r)\,:=\,r\A_n^\circ\cap\Z^{n^2},\end{equation}
it is seen that
\begin{equation}\label{ASMO}\ba{l}\!\!\!\ASMO(n,r):=\\[2mm]
\hs{1}\left\{\left(\ba{c@{\;}c@{\;}c}
A_{11}&\ldots&A_{1n}\\\vdots&&\vdots\\A_{n1}&\ldots&A_{nn}\ea\right)\!\!\in\Z^{n\times n}\:\left|\;\,\ba{l}
\bullet\ \sum_{j'=1}^n\!A_{ij'}=\sum_{i'=1}^n\!A_{i'\!j}=r\mbox{ \ for all }i,j\in[n]\!\\[2.2mm]
\bullet\ \sum_{j'=1}^jA_{ij'\!}\ge1\mbox{ \ for all }i\in[n],\;j\in[n\mi1]\\[2.2mm]
\bullet\ \sum_{j'=j}^nA_{ij'\!}\ge1\mbox{ \ for all }i\in[n],\;j\in[2,n]\\[2.2mm]
\bullet\ \sum_{i'=1}^iA_{i'\!j}\ge1\mbox{ \ for all }i\in[n\mi1],\;j\in[n]\\[2.2mm]
\bullet\ \sum_{i'=i}^nA_{i'\!j}\ge1\mbox{ \ for all }i\in[2,n],\;j\in[n]\ea\right.\right\}\\[20mm]
=\left\{\left(\ba{c@{\;}c@{\;}c}
A_{11}&\ldots&A_{1n}\\\vdots&&\vdots\\A_{n1}&\ldots&A_{nn}\ea\right)\!\!\in\Z^{n\times n}\:\left|\;\ba{l}
\bullet\ \sum_{j'=1}^n\!A_{ij'}=\sum_{i'=1}^n\!A_{i'\!j}=r\mbox{ \ for all }i,j\in[n]\\[2.2mm]
\bullet\ 1\le\sum_{j'=1}^j\!A_{ij'\!}\le r\mi1\mbox{ \ for all }i\in[n],\;j\in[n\mi1]\\[2.2mm]
\bullet\ 1\le\sum_{i'=1}^i\!A_{i'\!j}\le r\mi1\mbox{ \ for all }i\in[n\mi1],\;j\in[n]\ea\right.\!\right\}\!.\ea\!\!\!\!\!\!\end{equation}

Thus, $\ASMO(n,r)$ is the set of $n\t n$ integer-entry matrices
for which all complete row and column sums are $r$,
and all partial row and column sums extending from each end of the row or column are positive.
It follows that
\begin{equation}\label{ASMOASM}\ba{l}
\left\{\left(\ba{c@{\;\,}c@{\;\,}c}A_{11}\pl1&\ldots&A_{1n}\pl1\\\vdots&&\vdots\\A_{n1}\pl1&\ldots&A_{nn}\pl1\ea\right)
\;\left|\;\left(\ba{c@{\;}c@{\;}c}
A_{11}&\ldots&A_{1n}\\\vdots&&\vdots\\A_{n1}&\ldots&A_{nn}\ea\right)\in\ASM(n,r\mi n)\right\}\right.\subset\,\ASMO(n,r),\\[8mm]
\hs{85}\mbox{ for each }n,\,r\in\P\mbox{ with }r\ge n,\!\!\ea
\end{equation}
and that
\begin{equation}\label{ASMOempty}\ASMO(n,r)\,=\,\emptyset,\mbox{ \ for each }n,\,r\in\P\mbox{ with }r<n.\end{equation}
It can be seen immediately that the containment of~(\ref{ASMOASM}) is in fact an equality for $n=1$ and $n=2$, and
it follows from~(\ref{ASM3Alt}) in Section~8, that this is also the case for $n=3$.

We now state and prove the main result for the enumeration of $\ASM(n,r)$ and $\ASMO(n,r)$.
\newpage
\begin{theorem}For fixed $n\in\P$, there exists a function $\La(r)$, the Ehrhart polynomial
of the alternating sign matrix polytope $\A_n$, which satisfies:
\vs{-4}
\begin{enumerate}
\item $\La(r)$ is a polynomial in $r$ of degree $(n\mi1)^2$.
\item $|\ASM(n,r)|\,=\,\La(r)$, \ for each $r\in\N$.
\item $|\ASMO(n,r)|\,=\,(-1)^{n+1}\,\La(-r)$, \ for each $r\in\N$.
\item $\La(-1)\,=\,\La(-2)\,=\,\ldots\,=\La(-n\pl1)\,=\,0$.
\item $\ds\La(1)\,=\,\prod_{i=0}^{n\mi1}\!\frac{(3i\pl1)!}{(n\pl i)!}$.\vs{-4}
\item $\La(r)$ can be expressed in the form \ $\ds\La(r)\,=
\sum_{k=n-1}^{\;(n-1)^2}c_k\left(\,\ba{c}r+k\\[1.2mm](n\mi1)^2\ea\right)$,\\[3mm]
with $c_k\in\N$ for each $k\in[n\mi1,(n\mi1)^2]$.
\end{enumerate}
\end{theorem}
\textit{Proof.} \ All of the conclusions of this theorem follow straightforwardly from
results already obtained or stated in this paper.  First, the existence of an Ehrhart polynomial
$\mathcal{L}_{\!\A_n\!}(r)=\La(r)$ follows from the fact, implied by
Theorem~1, that~$\A_n$ is an integral convex polytope.  Properties (i)--(iii) in this theorem then follow
from (i)--(iii) of~(\ref{Ehr}) applied to $\A_n$, using~(\ref{dil}),~(\ref{dimA}) and~(\ref{ASMOdil}).
Property~(iv) follows from property~(iii) applied to~(\ref{ASMOempty}).
Property~(v) follows from property~(ii) applied to~(\ref{ASM1}).
Finally, property~(vi) follows
from~(iv) of~(\ref{Ehr}), and from property~(iv), which gives $c_0=c_1=\ldots=c_{n-2}=0$.\hspace{\fill}$\Box$

It follows that the explicit polynomial $\La(r)$ for a particular $n\in\P$ can be found by interpolation using
the $n\pl1$ values provided by $\La(0)=1$ and properties~(iv) and~(v), together with
$n^2\mi3n\pl1$ further values obtained by the direct enumeration of
cases of $\ASM(n,r)$ or $\ASMO(n,r)$ and the application of properties~(ii) and~(iii).
We have done this for $n=3$,~$4$ and~$5$, some of the required values being provided in Table~\ref{ASMnr}.
Together with the trivial cases $n=1$ and $2$, the results, expressed in the form of~(vi) of Theorem~2, are\ru{1.5}
\begin{equation}\label{ASM123}
A_1(r)=\biggl(\,\ba{c}r\\[0.3mm]0\ea\,\biggr),\quad A_2(r)=\biggl(\ba{c}r\pl1\\[0.3mm]1\ea\biggr),\quad
A_3(r)=\biggl(\ba{c}r\pl2\\[0.3mm]4\ea\biggr)+2\biggl(\ba{c}r\pl3\\[0.3mm]4\ea\biggr)+\biggl(\ba{c}r\pl4\\[0.3mm]4\ea\biggr),\ru{3.5}
\end{equation}
\begin{equation}\ba{l}
A_4(r)=\,3\biggl(\ba{c}r\pl3\\[0.3mm]9\ea\biggr)+80\biggl(\ba{c}r\pl4\\[0.3mm]9\ea\biggr)+415\biggl(\ba{c}r\pl5\\[0.3mm]9\ea\biggr)+
592\biggl(\ba{c}r\pl6\\[0.3mm]9\ea\biggr)+253\biggl(\ba{c}r\pl7\\[0.3mm]9\ea\biggr)+\mbox{}\\[4mm]
\qquad\qquad32\biggl(\ba{c}r\pl8\\[0.3mm]9\ea\biggr)+\biggl(\ba{c}r\pl9\\[0.3mm]9\ea\biggr),\ea\end{equation}
and
\begin{equation}\label{ASM5}\ba{l}
A_5(r)=\,70\biggl(\ba{c}r\pl4\\[0.3mm]16\ea\biggr)+14468\biggl(\ba{c}r\pl5\\[0.3mm]16\ea\biggr)+
521651\biggl(\ba{c}r\pl6\\[0.3mm]16\ea\biggr)+6002192\biggl(\ba{c}r\pl7\\[0.3mm]16\ea\biggr)+\mbox{}\\[4mm]
\qquad\qquad28233565\biggl(\ba{c}r\pl8\\[0.3mm]16\ea\biggr)+61083124\biggl(\ba{c}r\pl9\\[0.3mm]16\ea\biggr)+
64066830\biggl(\ba{c}r\pl10\\[0.3mm]16\ea\biggr)+\mbox{}\\[4mm]
\qquad\qquad32866092\biggl(\ba{c}r\pl11\\[0.3mm]16\ea\biggr)+
7998192\biggl(\ba{c}r\pl12\\[0.3mm]16\ea\biggr)+854464\biggl(\ba{c}r\pl13\\[0.3mm]16\ea\biggr)+\mbox{}\\[4mm]
\qquad\qquad34627\biggl(\ba{c}r\pl14\\[0.3mm]16\ea\biggr)+412\biggl(\ba{c}r\pl15\\[0.3mm]16\ea\biggr)+
\biggl(\,\ba{c}r\pl16\\[0.3mm]16\ea\,\biggr).\ea\end{equation}

We note that the fact that~(\ref{ASMOASM}) is an equality for $n=1$, 2 and~3, implies that
$|\ASMO(n,r)|=|\ASM(n,r\mi n)|$, for all $n\in\{1,2,3\}$ and $r\in\P$ with $r\ge n$, which in turn implies that
$A_n(r)=(-1)^{n+1}A_n(-n\mi r)$, for all $n\in\{1,2,3\}$ and $r\in\C$, as can be observed in~(\ref{ASM123}).

Finally, we now outline the results for $\B_n$ and $\SMS(n,r)$ which correspond
to those of this section for $\A_n$ and $\ASM(n,r)$.
Even though the results for $\B_n$ and $\SMS(n,r)$ are previously-known,
having been conjectured in~\cite{AnaDumGup66} and
first proved in~\cite{Ehr73,Sta73}, we state them here
to show their similarity with those for $\A_n$ and $\ASM(n,r)$.
The relative interior $\B^\circ_n$ of the Birkhoff polytope $\B_n$ is the set of positive real-entry $n\t n$
matrices for which all complete row and column sums are~$1$, and
$\SMSO(n,r):=r\B_n^\circ\cap\Z^{n^2}$ is the set of positive integer-entry matrices for which
all complete row and column sums are~$r$.  It follows that if $\SMS(n,r)$ and $\SMSO(n,r)$ are
substituted for $\ASM(n,r)$ and $\ASMO(n,r)$ in~(\ref{ASMOASM}) and~(\ref{ASMOempty}),
then the equations still hold, and furthermore that the containment in~(\ref{ASMOASM}) can
be replaced by an equality in all cases. Since $\B_n$ is an integral convex polytope,
it has an Ehrhart polynomial $\mathcal{L}_{\!\B_n\!}(r)$, which we denote~$H_n(r)$.
Applying~(\ref{Ehr}), and using~(\ref{SMS1}),~(\ref{dil}) and~(\ref{dimA}), it now follows that
all of Theorem~2 still holds if the substitutions~$H_n(r)$ for~$A_n(r)$, `Birkhoff polytope'
for `alternating sign matrix polytope', $\B_n$ for $\A_n$, $\SMS(n,r)$ for $\ASM(n,r)$,
$\SMSO(n,r)$ for $\ASMO(n,r)$, and $n!$ for $\prod_{i=0}^{n\mi1}(3i\pl1)!/(n\pl i)!$ are made.
Furthermore, the fact that the counterpart to~(\ref{ASMOASM}) is now an equality implies
that $|\SMSO(n,r)|=|\SMS(n,r\mi n)|$, for all $n,r\in\P$ with $r\ge n$, which leads to the further
properties that $H_n(r)=(-1)^{n+1}H_n(-n\mi r)$ for all $n\in\P$ and $r\in\C$, and, in the counterpart to~(vi) of Theorem~2, that
$c_k=c_{n(n-1)-k}$ for each $n\in\P$ and $k\in[n\mi1,(n\mi1)^2]$.
\font\tenmsb=msbm10 scaled \magstep1
\font\sevenmsb=msbm7 scaled \magstep1
\font\fivemsb=msbm5 scaled \magstep1
\newfam\msbfam
\textfont\msbfam=\tenmsb
\scriptfont\msbfam=\sevenmsb
\scriptscriptfont\msbfam=\fivemsb

\subsubsection{Higher Spin Alternating Sign Matrices of Size 3}
In the previous section, the enumeration of higher spin alternating sign matrices was studied using
a general, but nondirect, approach.
In this section, we consider the special case of $3\t3$ higher spin alternating sign matrices, and
provide a direct bijective derivation of the enumeration formula, from Theorem~2 and~(\ref{ASM123}),
\begin{equation}\label{ASM3}|\ASM(3,r)|\,=\biggl(\ba{c}r\pl2\\[0.3mm]4\ea\biggr)+
2\biggl(\ba{c}r\pl3\\[0.3mm]4\ea\biggr)+\biggl(\ba{c}r\pl4\\[0.3mm]4\ea\biggr),\mbox{ \ for each }r\in\N.\end{equation}

We begin by outlining a derivation of the corresponding semimagic squares formula,
\begin{equation}\label{SMS3}|\SMS(3,r)|\,=\biggl(\ba{c}r\pl2\\[0.3mm]4\ea\biggr)+
\biggl(\ba{c}r\pl3\\[0.3mm]4\ea\biggr)+\biggl(\ba{c}r\pl4\\[0.3mm]4\ea\biggr),\mbox{ \ for each }r\in\N,\end{equation}
since the derivation of~(\ref{ASM3}) will be similar.
The polynomial~(\ref{SMS3}) for $|\SMS(3,r)|$ was obtained directly in~\cite{AbrMos73,AnaDumGup66,Bon97,Mac15,NatIye72}.
The derivation given here most closely follows that of~\cite{Bon97}.

Let $S_3$ be the set of permutations of 123, and
consider the set
\begin{equation}\label{C}\ba{l}C(r):=\\[1.5mm]
\hs{6}\left\{\ba{l}a=(a_{123},a_{132},a_{213},\\[0.6mm]\hs{14}a_{231},a_{312},a_{321})\in\N^6\ea\:\left|\:\;\ba{l@{\;\:}l}
\bullet&\sum_{\s\in S_3}a_\s=r\\[2.2mm]
\bullet&a_{132}=0\mbox{ \ or \ }(a_{321}=0\mbox{ and }a_{132}\ne0)\\[1mm]
&\mbox{or \ }(a_{213}=0\mbox{ and }a_{132},\,a_{321}\ne0)
\ea\right.\right\}.\ea\end{equation}
Using the fact that, for any $m,n\in\N$, the number of tuples $(a_1,\ldots,a_{m+n})\in\P^m\t\N^n$ with
$\sum_{i=1}^{m+n}a_i=r$ is $\biggl(\ba{c}r\pl n\mi1\\[0.1mm]m\pl n\mi1\ea\biggr)$ (which can be obtained by considering
tuples $(b_1,\ldots,b_{m+n-1})\in\P^{m+n-1}$ with $b_1<\!\ldots\!<b_m\le b_{m+1}\le\!\ldots\!\le b_{m+n-1}\le r$, and
using the bijection $b_i=\sum_{j=1}^ia_j$ for each $i\in[m\pl n\mi1]$), it follows that
\begin{equation}\label{NC}|C(r)|\,=\biggl(\ba{c}r\pl4\\[0.3mm]4\ea\biggr)+
\biggl(\ba{c}r\pl3\\[0.3mm]4\ea\biggr)+\biggl(\ba{c}r\pl2\\[0.3mm]4\ea\biggr),\end{equation}
where the three binomial coefficients correspond respectively to the three alternatives in the second
condition of~(\ref{C}).

Now let $P_\s$ be the permutation matrix given
by $(P_\s)_{ij}:=\d_{i,\s_j}$ for each $\s\in S_3$ and $i,j\in[3]$, and define, for each $a\in C(r)$,
\begin{equation}\label{theta}\theta(a)\,:=\ds\sum_{\s\in S_3}\!a_\s P_\s\,=\,\left(\ba{c@{\;\,}c@{\;\,}c}
a_{123}\pl a_{132}&a_{213}\pl a_{312}&a_{231}\pl a_{321}\\
a_{213}\pl a_{231}&a_{123}\pl a_{321}&a_{132}\pl a_{312}\\
a_{312}\pl a_{321}&a_{132}\pl a_{231}&a_{123}\pl a_{213}\ea\right).\end{equation}
It can immediately be seen that $\theta(a)\in\SMS(3,r)$ for each $a\in C(r)$.
It can also be checked straightforwardly that $\theta$
is a bijection between $C(r)$ and $\SMS(3,r)$, the inverse mapping being, for each $A\in\SMS(3,r)$,
\begin{equation}\label{thetainv}\theta^{-1}(A)\,=\left\{\ba{l}
(A_{11},0,A_{33}\mi A_{11},A_{32},A_{23},A_{22}\mi A_{11}),\ \ A_{11}\le A_{22}\mbox{ and }A_{11}\le A_{33}\\[2.1mm]
(A_{22},A_{11}\mi A_{22},A_{33}\mi A_{22},A_{13},A_{31},0),\ \ A_{22}<A_{11}\mbox{ and }A_{22}\le A_{33}\\[2.1mm]
(A_{33},A_{11}\mi A_{33},0,A_{21},A_{12},A_{22}\mi A_{33}),\ \ A_{33}<A_{11}\mbox{ and }A_{33}<A_{22},\ea\right.\!\!\end{equation}
where the three cases of~(\ref{thetainv}) correspond respectively to the three alternatives in the second
condition of~(\ref{C}).
Therefore $|\SMS(3,r)|=|C(r)|$, so that~(\ref{NC}) gives~(\ref{SMS3}), thereby completing its bijective derivation.

Proceeding to $3\times3$ higher spin alternating sign matrices,
we first note that $\ASM(3,r)$ is simply the set of $3\times 3$ integer-entry matrices
for which all complete row and column sums are $r$,
and the eight external entries are nonnegative,
\begin{equation}\label{ASM3Alt}\ba{l}\!\!\ASM(3,r)\,=\\[2mm]
\hs{2.9}\left\{\left.\!A\!=\!\!\left(\ba{c@{\;}c@{\;}c}
A_{11}&A_{12}&A_{13}\\
A_{21}&A_{22}&A_{23}\\
A_{31}&A_{32}&A_{33}\ea\right)\!\in\Z^{3\times 3}\:\right|\;\ba{l}
\bullet\ \sum_{j'=1}^3\!A_{ij'}=\sum_{i'=1}^3\!A_{i'\!j}=r\mbox{ \ for all }i,j\in[3]\\[2.5mm]
\bullet\ A_{ij}\ge0\mbox{ \ for all }(i,j)\in[3]\t[3]/\{(2,2)\}\ea\right\}\!.\ea\!\!\end{equation}
To obtain (\ref{ASM3Alt}) from (\ref{ASM}), it only needs to be
verified that if $A$ is an element of the RHS of (\ref{ASM3Alt}),
then each partial sum $A_{12}+A_{22}$, $A_{21}+A_{22}$, $A_{23}+A_{22}$
and $A_{32}+A_{22}$ is nonnegative. For the case of $A_{12}+A_{22}$,
$A_{12}\pl A_{22}\pl A_{32}=A_{31}\pl A_{32}\pl A_{33}$ $(=r)$
gives $A_{12}\pl A_{22}=A_{31}\pl A_{33}$, but $A_{31}\ge0$
and $A_{33}\ge0$ now imply that $A_{12}\pl A_{22}\ge0$ as required.
The other three cases follow by symmetry.

The set of standard alternating sign matrices, $\ASM(3,1)$, is the set $\{P_\s\mid\s\in S_3\}$ of $3\t3$ permutation
matrices, together with the matrix\rule[-3.5ex]{0ex}{8.5ex} $\left(\mbox{\footnotesize$\ba{r@{\;}r@{\;\;\,}r}0&1&0\\1&-1&1\\0&1&0\ea$}\right)$ which
will be denoted as $P_0$. Setting $S'_3:=S_3\cup\{0\}$, we now define, in analogy with~(\ref{C}) and~(\ref{theta}),
\begin{equation}\label{CC}\ba{l}C'(r):=\\[1.5mm]
\hs{5}\left\{\ba{l}a=(a_{123},a_{132},a_{213},a_{231},\\[0.6mm]\hs{15.8}a_{312},a_{321},a_0)\in\N^7\ea\:\left|\:\;\ba{l@{\;\:}l}
\bullet&\sum_{\s\in S'_3}a_\s=r\\[2.2mm]
\bullet&a_{132}=a_0=0 \mbox{ \ or}\\[1mm]
&(a_{321}=a_0=0\mbox{ and }a_{132}\ne0)\mbox{ \ or}\\[1mm]
&(a_{213}=a_0=0\mbox{ and }a_{132},\,a_{321}\ne0)\mbox{ \ or}\\[1mm]
&(a_{123}=a_{321}=0\mbox{ and }a_0\ne0)
\ea\right.\right\},\ea\end{equation}
and, for each $a\in C'(r)$,
\begin{equation}\theta'(a)\,:=\ds\sum_{\s\in S'_3}\!a_\s P_\s\,=\,\left(\ba{c@{\;\;\:}c@{\;\;\:}c}
a_{123}\pl a_{132}&a_{213}\pl a_{312}\pl a_0&a_{231}\pl a_{321}\\
a_{213}\pl a_{231}\pl a_0&a_{123}\pl a_{321}\mi a_0&a_{132}\pl a_{312}\pl a_0\\
a_{312}\pl a_{321}&a_{132}\pl a_{231}\pl a_0&a_{123}\pl a_{213}\ea\right).\end{equation}
It immediately follows that
\begin{equation}\label{NCC}|C'(r)|\,=\biggl(\ba{c}r\pl4\\[0.3mm]4\ea\biggr)+
2\biggl(\ba{c}r\pl3\\[0.3mm]4\ea\biggr)+\biggl(\ba{c}r\pl2\\[0.3mm]4\ea\biggr),\end{equation}
and that $\theta'(a)\in\ASM(3,r)$ for each $a\in C'(r)$.
It can also be checked, using~(\ref{ASM3Alt}) and the bijection~(\ref{theta}) and (\ref{thetainv}) between
$C(r)$ and $\SMS(3,r)$, that $\theta'$
is a bijection between $C'(r)$ and $\ASM(3,r)$, the inverse mapping being, for each $A\in\ASM(3,r)$,
\begin{equation}\label{thetadinv}\ba{l}\theta'^{\,-1}(A)=\\[2mm]
\hs{6}\left\{\ba{l}
(A_{11},0,A_{33}\mi A_{11},A_{32},A_{23},A_{22}\mi A_{11},0),\ \ A_{11}\le A_{22}\mbox{ and }A_{11}\le A_{33}\\[2.1mm]
(A_{22},A_{11}\mi A_{22},A_{33}\mi A_{22},A_{13},A_{31},0,0),\ \ 0\le A_{22}<A_{11}\mbox{ and }A_{22}\le A_{33}\\[2.1mm]
(A_{33},A_{11}\mi A_{33},0,A_{21},A_{12},A_{22}\mi A_{33},0),\ \ A_{33}<A_{11}\mbox{ and }A_{33}<A_{22}\\[2.1mm]
(0,A_{11},A_{33},A_{13},A_{31},0,-A_{22}),\ \ A_{22}<0,\ea\right.\ea\end{equation}
where the four cases of~(\ref{thetadinv}) correspond respectively to the four alternatives in the second
condition of~(\ref{CC}).
Therefore $|\ASM(3,r)|=|C'(r)|$, so that~(\ref{NCC}) gives~(\ref{ASM3}), thereby completing its bijective derivation.

In this section, particular sets $C(r)$ and $C'(r)$ of weak compositions of $r$
into six and seven parts, with at most five parts nonzero,
have been found which are in bijection with $\SMS(3,r)$ and $\ASM(3,r)$ respectively,
where the mapping from a weak composition to a matrix simply uses the parts of the composition
as coefficients in a linear combination of permutation matrices or standard alternating sign matrices.
It would clearly be worth investigating whether these results can be generalized to larger integers $n$ to give
sets of weak compositions of $r$ into~$n!$ and $\prod_{i=0}^{n\mi1}(3i\pl1)!/(n\pl i)!$
parts, with at most $(n\mi1)^2\pl1$ parts nonzero, which are similarly in bijection with $\SMS(n,r)$ and $\ASM(n,r)$
respectively.  This, however, seems to be a challenging problem.

\subsubsection{Discussion}
We conclude this paper by outlining some possible directions for further research.
While all of these can be simply stated, and follow from the content of this paper in obvious ways,
some may well turn out to be of limited interest, or to lead to prohibitively-difficult problems.

In this paper, we studied the configurations of certain integrable statistical mechanical models with
domain wall boundary conditions.
One area for additional
work would therefore be the consideration of further integrable statistical mechanical models with
boundary conditions of domain-wall-type.
A possible component of such work would be the derivation of determinant formulae for the partition functions of such models,
this having already been done for certain models in~\cite{CarFodWheZup07,DowFod06}.

In Section~3, we discussed the determinant formula for
higher spin vertex models.  Although it seems unlikely that this formula can be used for the unweighted enumeration of
higher spin alternating sign matrices or semimagic squares,
it may still be possible for it to be used to obtain certain weighted enumeration formulae.  More specifically,
each vertex type of $\V(r)$, as defined in~(\ref{Vr}), would be given a weight, and the overall weight of a higher spin
alternating sign matrix would be the product of the weights of all the vertex types associated with the corresponding edge matrix pair.
Certain such weighted enumeration formulae  are already
known for standard alternating sign matrices (see for example~\cite{ColPro05,ElkKupLarPro92,Kup96}). In these cases,
vertex type~(4) in Figure~\ref{19V}
(which corresponds to each entry of~$-1$ in a standard alternating sign matrix) is given a weight of~$2$ or~$3$, and the other five
vertex types of $\V(1)$ are each given a weight of~1.

In Sections~4 and~5, we considered various other combinatorial objects, namely certain sets of lattice paths, monotone
triangles, corner sum matrices and complementary edge matrix pairs, which are in simple bijection with higher spin
alternating sign matrices.  Each of these representations could be studied further, and
might provide useful statistics according to which higher spin alternating sign matrices
could be weighted and classified.  The fully packed loop configurations defined in Section~5 also
seem worthy of further study, especially since their classification according to the
link pattern formed among the external points by the open paths, may provide an interesting generalization of
such a classification for the standard case of $r=1$ (see for example~\cite{CasKra04,CasKraLasNad05,DifZinZub04,DifZub04,Wie00,Zub04}).

It is also common to impose further conditions on semimagic squares (see for example~\cite{BecCohCuoGri03} and \cite[Proposition~4.6.21]{Sta86})
and standard alternating sign matrices (see for example~\cite{Kup96,Oka06,RazStr04,Rob91,Rob00}).
It would therefore seem natural to investigate the effects of imposing similar such conditions on higher spin alternating
sign matrices, it being expected that this would lead to the introduction of further, not necessarily integral, convex
polytopes.

Semimagic squares are a special case of contingency or frequency tables (see for example~\cite{GaiMan77}),
these simply being nonnegative integer-entry rectangular matrices with arbitrary prescribed row and column sums.
Similar generalizations of higher spin alternating sign matrices in which the complete row and column sums are prescribed, but not necessarily
all equal, and the matrices are not necessarily square, would thus provide another direction for further study.

In Section~6, we defined the alternating sign matrix polytope, and showed that its vertices
are standard alternating sign matrices.  Many additional results for this polytope, including
the enumeration of its facets, the characterization
of its face lattice, and the description of its projection to the permutohedron have been obtained in~\cite{Str07},
but this polytope could still be investigated further.

In Section~7, we obtained certain enumerative formulae for higher spin alternating sign matrices with fixed size $n$ and
variable line-sum $r$, but
we have not obtained formulae for the case of variable $n$ and fixed~$r$.
Formulae for
semimagic squares with variable $n$, and $r=2$ or $r=3$ are known
(see for example~\cite{AnaDumGup66},~\cite[Ex.~25, Page~124]{Com74} and~\cite[Sec.~5.5]{Sta99}, but note that
the $r=3$ formula in~\cite{Com74} contains misprints).  However, it might be expected that the
difference in difficulties in obtaining, say, the $r=2$ semimagic square formula and an $r=2$ higher spin
alternating sign matrix formula, might be comparable to the substantial difference in difficulties in obtaining the $r=1$
formulae~(\ref{SMS1}) and~(\ref{ASM1}).

Finally, as indicated at the end of Section~8, it would be interesting to see whether the approach used there
can be generalized to give bijective derivations of enumerative formulae for semimagic squares or
higher spin alternating sign matrices with fixed size larger
than~3.\\[10mm]
\textbf{Acknowledgements}\\[6mm]
We thank Omar Foda and Jessica Striker for correspondence regarding
their papers~\cite{CarFodKit06} and~\cite{Str07}.
VK is supported by an EPSRC-funded postgraduate scholarship.

\end{document}